\documentclass[twoside,11pt,reqno]{amsart}
\usepackage{amsmath,amssymb,amscd,mathrsfs}

\makeatletter

\@addtoreset{equation}{section}

\newtheorem{Proposition}{Proposition}[section]
\newtheorem{Lemma}[Proposition]{Lemma}
\newtheorem{Theorem}[Proposition]{Theorem}
\newtheorem{Corollary}[Proposition]{Corollary}

\newtheorem{Remark}[Proposition]{Remark}

\newtheorem{Example}[Proposition]{Example}

\newbox\squ  
\setbox\squ=\hbox{\vrule width.6pt
   \vbox{\hrule height.6pt width.4em\kern1ex
         \hrule height.6pt}%
   \vrule width.6pt}

\def\op{\operatorname{op}}

\def\V{\mathbb{V}}
\def\G{\mathbb{G}}
\def\ef{\mathbb{F}}
\def\ee{\mathbb{E}}

\def\sgn{\operatorname{sgn}}

\def\Mod{\operatorname{-mod}}

\def\id{\operatorname{id}}

\def\Id{\operatorname{Id}}
\def\C{{\mathbb C}}

\def\Z{{\mathbb Z}}

\def\0{{\bar 0}}
\def\1{{\bar 1}}

\def\hom{{\operatorname{Hom}}}

\def\End{{\operatorname{End}}}
\def\Adj{{\operatorname{Adj}}}

\def\bid{\hbox{\boldmath{$1$}}}
\def\eps{{\varepsilon}}
\def\phi{{\varphi}}

\def\underbar{\mathpalette\@underbar}
\def\@underbar#1#2{\settowidth{\@tempdimb}{$#1#2$}\@tempdimb=0.8\@tempdimb
                   \ooalign{$#1#2$\crcr%
                         \hfil\rule[-.5mm]{\@tempdimb}{.4pt}\hfil}}

\newdimen\hoogte    \hoogte=8pt    
\newdimen\breedte   \breedte=8pt   
\newdimen\dikte     \dikte=0.5pt    

\newenvironment{young}{\begingroup
       \def\vr{\vrule height0.8\hoogte width\dikte depth 0.2\hoogte}
       \def\fbox##1{\vbox{\offinterlineskip
                    \hrule height\dikte
                    \hbox to \breedte{\vr\hfill##1\hfill\vr}
                    \hrule height\dikte}}
       \vbox\bgroup \offinterlineskip \tabskip=-\dikte \lineskip=-\dikte
            \halign\bgroup &\fbox{##\unskip}\unskip  \crcr }
       {\egroup\egroup\endgroup}
\def\diagram#1{\relax\ifmmode\vcenter{\,\begin{young}#1\end{young}\,}\else%
              $\vcenter{\,\begin{young}#1\end{young}\,}$\fi}

\newdimen\Hoogte    \Hoogte=12pt    
\newdimen\Breedte   \Breedte=12pt   
\newdimen\Dikte     \Dikte=0.5pt    

\newenvironment{Young}{\begingroup
       \def\vr{\vrule height0.8\Hoogte width\Dikte depth 0.2\Hoogte}
       \def\fbox##1{\vbox{\offinterlineskip
                    \hrule height\Dikte
                    \hbox to \Breedte{\vr\hfill##1\hfill\vr}
                    \hrule height\Dikte}}
       \vbox\bgroup \offinterlineskip \tabskip=-\Dikte \lineskip=-\Dikte
            \halign\bgroup &\fbox{##\unskip}\unskip  \crcr }
       {\egroup\egroup\endgroup}
\def\Diagram#1{\relax\ifmmode\vcenter{\,\begin{Young}#1\end{Young}\,}\else%
              $\vcenter{\,\begin{Young}#1\end{Young}\,}$\fi}

\begin{document}

\title[Springer fibers]{\sc Symmetric functions, Parabolic Category $\mathcal O$ and the Springer Fiber}
\author{\sc Jonathan Brundan}
\begin{abstract}
We prove that the center of a regular block of parabolic category
$\mathcal O$ for the general linear Lie algebra
is isomorphic to the cohomology algebra of a corresponding
Springer fiber. This was conjectured by Khovanov.
We also find presentations for the centers of singular blocks,
which are cohomology algebras of Spaltenstein varieties.
\end{abstract}
\thanks{Research supported in part by NSF grant no. DMS-0139019.}
\address{Department of Mathematics, University of Oregon, Eugene, Oregon, USA.}
\email{brundan@uoregon.edu}

\maketitle

\section{Introduction}

Fix a natural number $n$.
By a {\em composition}, respectively, a {\em partition} of $n$
we mean a tuple $\lambda = (\lambda_i)_{i \in\Z}$,
respectively, a sequence
$\lambda = (\lambda_1 \geq \lambda_2 \geq \cdots)$
of non-negative integers summing to $n$.
For a composition $\nu$ of $n$, we write 
$$
S_\nu = \cdots \times S_{\nu_1} \times S_{\nu_2} \times \cdots
$$ 
for the usual
parabolic subgroup of the symmetric group $S_n$ parametrized by $\nu$,
and call $\nu$ {\em regular} if $S_\nu = \{1\}$.
Let $P :=\C[x_1,\dots,x_n]$,
viewed as a graded commutative algebra with each $x_i$ in degree two.
The symmetric group $S_n$ acts 
as usual on $P$ by permuting the generators.
Let $P_\nu$ be the subalgebra 
$\C[x_1,\dots,x_n]^{S_\nu}$ 
of all
$S_\nu$-invariants in $P$.
Given distinct integers $i_1,\dots,i_m$, we 
write $e_r(\nu;i_1,\dots,i_m)$ and $h_r(\nu;i_1,\dots,i_m)$ for the $r$th
elementary and complete symmetric polynomials in variables
$X_{i_1}\cup\cdots\cup X_{i_m}$, where
$$
X_i := \left\{x_j\:\Bigg|\: \sum_{h < i} \nu_h < j \leq \sum_{h \leq i} \nu_h\right\}.
$$
Note $P_\nu$ is  itself a free polynomial algebra of rank $n$
generated by the elements $\{e_r(\nu;i)\mid i \in \Z, 1 \leq r \leq \nu_i\}$
and also by
$\{h_r(\nu;i)\mid i \in\Z, 1 \leq r \leq \nu_i\}$.
Moreover, 
we have that 
$e_r(\nu;i_1,\dots,i_m) = 0$ for $r > \nu_{i_1}+\cdots+\nu_{i_m}$.

Now fix a composition $\mu$ of $n$ and
let $\lambda$ denote the
{\em transpose partition}, 
that is, $\lambda_j$ counts the number of $i \in \Z$
such that $\mu_i \geq j$.
Let $I^\mu_\nu$ denote the 
homogeneous ideal of $P_\nu$ generated by
\begin{equation*}
\left\{h_r(\nu;i_1,\dots,i_m)\:\bigg|
\begin{array}{l}
m \geq 1, \:\:i_1,\dots,i_m\text{ distinct integers},\\
r > \lambda_1+\cdots+\lambda_m - \nu_{i_1}-\cdots-\nu_{i_m}\end{array}\right\}.
\end{equation*}
Equivalently (see Lemma~\ref{l2}), $I^\mu_\nu$ is the ideal generated by
\begin{equation*}
\left\{e_r(\nu;i_1,\dots,i_m)\:\Bigg|
\begin{array}{l}
m \geq 1, \:\: i_1,\dots,i_m\text{ distinct integers},\\
r > \nu_{i_1}+\cdots+\nu_{i_m} - \lambda_{l+1}-\lambda_{l+2}-\cdots\\
\text{where }l:=\#\{i \in \Z\mid \nu_i > 0, i \neq i_1,\dots,i_m\}
\end{array}\right\}.
\end{equation*}
Let $C^\mu_\nu$ 
denote the graded quotient $P_\nu / I^\mu_\nu$.
These algebras have 
natural geometric realizations, as follows.
\begin{itemize}
\item
If both $\mu$ and $\nu$ are regular
then $I^\mu_\nu$ is simply the ideal of $P$
generated by all homogeneous symmetric polynomials of positive degree.
In this case, we write simply $I$ and $C$ for $I^\mu_\nu$ and $C^\mu_\nu$.
The algebra $C$ is the {\em coinvariant algebra},
which by a classical theorem of Borel
is isomorphic to the cohomology algebra 
 (with complex coefficients)
of the {\em flag manifold} $F$ of complete flags in $\C^n$.
\item
If just $\mu$ is regular, we denote $I^\mu_\nu$ and $C^\mu_\nu$ simply
by $I_\nu$ and $C_\nu$.
The algebra $C_\nu$ is isomorphic to 
the subalgebra $C^{S_\nu}$ of all $S_\nu$-invariants in $C$,
which is the cohomology algebra of the {\em partial flag manifold} $F_\nu$
of flags $\cdots \subseteq V_1 \subseteq V_{2} \subseteq \cdots$
with $\dim V_j = \sum_{i \leq j} \nu_i$ for each $j \in \Z$.
\item
If just $\nu$ is regular, we denote $I^\mu_\nu$ and $C^\mu_\nu$ simply
by $I^\mu$ and $C^\mu$.
The ideal $I^\mu$ is generated by 
the elementary symmetric functions $e_r(x_{i_1},\dots,x_{i_m})$ for every $m \geq 1$,
$1 \leq i_1 < \cdots < i_m \leq n$ and $r > m-\lambda_{n-m+1}-\lambda_{n-m+2}-\cdots$.
This is Tanisaki's presentation \cite{T} (simplifying De Concini and Procesi's 
original
work \cite{DP})
for the cohomology algebra of the {\em Springer fiber}
$F^\mu$ of all flags in $F$
stabilized by the nilpotent matrix $x_\mu$ of Jordan type $\mu$.
\item
Generalizing these special cases, we will prove in 
a subsequent article \cite{BO} (see also Remark~\ref{drem} 
and Example~\ref{oneex} below) that 
the algebra $C^\mu_\nu$ 
for arbitrary $\mu$ and $\nu$
is isomorphic to the cohomology algebra
of the {\em Spaltenstein variety} $F^\mu_\nu$ introduced in \cite{Sp},
that is, the subvariety of $F_\nu$ consisting 
of all partial flags
$\cdots \subseteq V_1 \subseteq V_2 \subseteq \cdots$ such that
$x_\mu V_j \subseteq V_{j-1}$ for each $j$.
\end{itemize}

The main goal of this article
is to explain how the algebras $C_\nu^\mu$ arise in representation 
theory. Consider 
the Lie algebra $\mathfrak{g} = \mathfrak{gl}_n(\C)$.
Let $\mathfrak{h}$ be the
Cartan subalgebra of diagonal matrices
and $\mathfrak{b}$ be the Borel subalgebra of upper triangular matrices.
Let $\eps_1,\dots,\eps_n$ be the basis for the 
vector space $\mathfrak{h}^*$
that is dual to the standard basis
$e_{1,1},\dots,e_{n,n}$ of $\mathfrak{h}$ 
consisting of matrix units.
The Bernstein-Gelfand-Gelfand category $\mathcal O$ 
introduced originally in \cite{BGG}
is the category of all finitely generated
$\mathfrak{g}$-modules that are
locally finite over $\mathfrak{b}$ and semisimple over $\mathfrak{h}$.
For $\alpha \in \mathfrak{h}^*$, write $L(\alpha)$ for the
irreducible highest weight module of highest weight $(\alpha - \rho)$
where $\rho := -\eps_2-2\eps_3-\cdots-(n-1)\eps_n$\footnote{This choice of
$\rho$ is congruent to the usual choice (half the sum of the positive roots)
modulo $\eps_1+\eps_2+\cdots+\eps_n$.
The choice here is more convenient when working with
$\mathfrak{gl}_n(\C)$ rather than $\mathfrak{sl}_n(\C)$ because it is
itself an integral weight.}.
These are the irreducible objects in $\mathcal O$.

For $\nu$ as above,  let $\mathcal O_\nu$ be the Serre subcategory of
$\mathcal O$ generated by the irreducible modules $L(\alpha)$
for all $\alpha = \sum_{i=1}^n a_i \eps_i \in \mathfrak{h}^*$ such that exactly
$\nu_i$ of the coefficients $a_1,\dots,a_n$ are equal to $i$ for each $i \in \Z$.
This is an {\em integral block} of $\mathcal O$.
For $\mu$ as above, let $\mathfrak{p}$ be the standard parabolic subalgebra
of block upper triangular matrices with blocks of size
$\dots,\mu_{1},\mu_2,\dots$ down the diagonal.
Let $\mathcal O^\mu$ be the full subcategory of $\mathcal O$ consisting of all modules
that are locally finite over $\mathfrak{p}$. This is 
{\em parabolic category $\mathcal O$}.
Finally we set $\mathcal O^\mu_\nu := \mathcal O^\mu \cap \mathcal O_\nu$, 
an integral
block\footnote{The fact that
$\mathcal O^\mu_\nu$ really is an 
indecomposable subcategory of $\mathcal O^\mu$,
so is a block of $\mathcal O^\mu$ in the usual sense, is
explained in the discussion immediately following the statement of Theorem 2 in
the introduction of \cite{cyclo}.} of parabolic category $\mathcal O$.
Recall the center of an additive category $\mathcal C$ is the
ring $Z(\mathcal C)$ of all natural transformations from the identity functor
to itself. 

\vspace{1mm}\noindent{\bf Main Theorem.}
{\em $Z(\mathcal O^\mu_\nu) \cong C^\mu_\nu$.}
\vspace{1mm}

To make this isomorphism explicit, let $z_1,\dots,z_n$
be the generators for the center
$Z(\mathfrak{g})$
of the universal enveloping algebra $U(\mathfrak{g})$ of $\mathfrak{g}$
 determined by the property
that, for $\alpha = \sum_{i=1}^n a_i \eps_i \in \mathfrak{h}^*$,
the element $z_r$ acts on the irreducible module $L(\alpha)$
by the scalar $e_r(a_1,\dots,a_n)$, the $r$th elementary symmetric polynomial
evaluated at $a_1,\dots,a_n$.
Let 
$p_\nu:Z(\mathfrak{g}) \rightarrow Z(\mathcal O_\nu)$ denote the canonical homomorphism sending $z$ to the natural transformation defined by left 
multiplication by $z$.
There is also a homomorphism $q_\nu:Z(\mathfrak{g}) \rightarrow C_\nu$
with
$q_\nu(z_r) = e_r(x_1+a_1,\dots,x_n+a_n)$
for each $r =1,\dots,n$, where $a_1,\dots,a_n$ here are defined so
that $a_1 \leq \cdots \leq a_n$ and as before
exactly $\nu_i$ of the integers $a_1,\dots,a_n$
are equal to $i$ for each $i \in \Z$.
By a famous result of Soergel \cite{Soergel} 
(see $\S$6 below for precise references) both of the maps $p_\nu$ and $q_\nu$ are surjective,
and there is a unique isomorphism
$c_\nu$ making the following diagram commute:
$$
\begin{CD}
&\!Z(\mathfrak{g})\: \\
\\
Z(\mathcal O_\nu)
@>\sim >c_\nu> &C_\nu.
\end{CD}
\begin{picture}(0,0)
\put(-66,-1){\makebox(0,0){$\swarrow$}}
\put(-68,7){\makebox(0,0){$_{p_\nu}$}}
\put(-66,-1){\line(1,1){14}}
\put(-16,-1){\makebox(0,0){$\searrow$}}
\put(-14.5,7){\makebox(0,0){$_{q_\nu}$}}
\put(-16,-1){\line(-1,1){14}}
\end{picture}
$$
We actually prove that there exists a unique isomorphism
$c^\mu_\nu$
making the following diagram commute:
$$
\begin{CD}
Z(\mathcal O_\nu)&@>\sim>c_\nu>&C_\nu\\
@Vr^\mu_\nu VV&&@VVs^\mu_\nu V\\
Z(\mathcal O^\mu_\nu)&@>\sim>c^\mu_\nu>&C^\mu_\nu,
\end{CD}
$$
where $r^\mu_\nu:Z(\mathcal O_\nu) \rightarrow Z(\mathcal O^\mu_\nu)$ is the 
restriction map arising from the embedding of $\mathcal O^\mu_\nu$
into $\mathcal O_\nu$ and $s^\mu_\nu:C_\nu \twoheadrightarrow C^\mu_\nu$ 
is the canonical quotient
map, which exists because $I_\nu \subseteq I^\mu_\nu$.
Note the surjectivity of the map $r^\mu_\nu$, an essential step in the proof, 
was already established in \cite[Theorem 2]{cyclo}.

In the special case that $\nu$ is regular, this proves that
the center of a regular block of $\mathcal O^\mu$ is isomorphic
to the cohomology algebra of the Springer fiber $F^\mu$,
exactly as was conjectured by Khovanov in \cite[Conjecture 3]{Kh}\footnote{Of the other conjectures from Khovanov's paper, 
\cite[Conjecture 4]{Kh} was already proved
by Mazorchuk and Stroppel \cite[Theorem 5.2]{MS}, and \cite[Conjecture 2]{Kh} 
then follows by \cite[Theorem 2]{St}. This just leaves the geometric 
\cite[Conjecture 1]{Kh}, which appears still to be unresolved.}.
Using the surjectivity 
of the map $r^\mu_\nu$ established in \cite{cyclo}, Stroppel 
\cite[Theorem 1]{St} has also independently
found a proof of Khovanov's conjecture.
Stroppel's approach based on deformation
is quite different from the 
strategy followed 
here.
The approach here 
has the advantage of yielding 
at the same time the explicit description of the centers of all
singular blocks. By \cite[Theorem 5.11]{cyclo}, these can also be reinterpreted
in terms of the centers of blocks of 
degenerate cyclotomic Hecke algebras.

The key idea in the proof 
is the 
construction of an action of the general linear Lie algebra
$\hat{\mathfrak{g}} := \mathfrak{gl}_{\infty}(\C)$ on the 
direct sum $\bigoplus_\nu Z(\mathcal O_\nu) \cong \bigoplus_\nu C_\nu$
of the centers of all
integral blocks of $\mathcal O$.
In this construction, the Chevalley generators
of $\hat{\mathfrak{g}}$ 
act as the {trace maps} in the sense of \cite{Btr} associated to some
canonical adjunctions between the
 special translation functors that arise by tensoring with the natural $\mathfrak{g}$-module and its dual.
In a similar way, there is an action of
$\hat{\mathfrak{g}}$ on the direct sum
$\bigoplus_\nu Z(\mathcal O_\nu^\mu) \cong \bigoplus_\nu 
C^\mu_\nu$ of the centers of all integral blocks
of $\mathcal O^\mu$ such that
the 
canonical map
$\oplus_\nu r^\mu_\nu$
is a $\hat{\mathfrak{g}}$-module homomorphism
(see Theorem~\ref{dthm}).
As explained in more detail in \cite{BO}, this action is closely related to
Ginzburg's geometric construction of representations of
the general linear group \cite{Gi}; see also \cite{BG} and \cite[$\S$7]{Gi2}.

The remainder of the article is organized as follows.
In $\S$2, we begin with some preliminaries on symmetric functions,
then deduce some elementary properties of the algebras $C^\mu_\nu$.
In $\S$3, we construct an action of $\hat{\mathfrak{g}}$ on the direct sum
$\bigoplus_\nu C_\nu$ of all partial coinvariant algebras by 
exploiting Schur-Weyl duality, 
paralleling the idea of
Braverman and Gaitsgory \cite{BG} on the geometric side.
In $\S$4, we use this action to deduce the dimension of 
$C^\mu_\nu$ 
from known properties of $C^\mu$.
Then in $\S$5 we reinterpret
the actions of the Chevalley generators of 
$\hat{\mathfrak{g}}$
on $\bigoplus_\nu C_\nu$ as certain 
trace maps.
Only in $\S$6 do we finally relate things back to category 
$\mathcal O$, using the full strength of Soergel's theory from \cite{Soergel}
to complete the proof of the Main Theorem.

\vspace{2mm}
\noindent
{\em Acknowledgements.}
I would like to thank Catharina Stroppel, Victor Ostrik and Nick Proudfoot 
for many helpful 
discussions, and the referees for some helpful comments.

\section{Preliminaries}

Let $\lambda$ be a partition of $n$.
Recall that the {\em Young diagram} of $\lambda$
consists of $\lambda_i$ boxes in its $i$th row; 
for example, the Young diagram of $\lambda = (4\,3^2\,2)$ is
$$
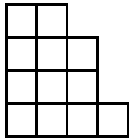
\Diagram{&\cr&&\cr&&\cr&&&\cr}
$$
Given in addition a composition $\nu$ of $n$,
a {\em $\lambda$-tableau} of {\em type $\nu$}
means some filling of the boxes of this diagram
by integers so that there are exactly $\nu_i$ entries equal to $i$
for each $i \in \Z$.
A $\lambda$-tableau is {\em column strict} if its entries are strictly
increasing in each column from bottom to top, and it is {\em standard}
if it is column strict and 
in addition its entries are weakly increasing in each row from left to right.
The {\em Kostka number} $K_{\lambda,\nu}$ is the number of standard $\lambda$-tableaux of type $\nu$. 
It is well known that $K_{\lambda,\nu}$ is non-zero if and only if
$\lambda \geq \nu^+$,
where $\geq$ denotes the usual dominance ordering on partitions
and $\nu^+$
is the unique partition of $n$ whose non-zero parts have the same
multiplicities as in the composition $\nu$;
see for example \cite[(I.6.5)]{Mac}.

Let $e_r(x_1,\dots,x_n)$ and $h_r(x_1,\dots,x_n)$ denote
the $r$th elementary and complete symmetric polynomials
in commuting variables $x_1,\dots,x_n$,
adopting the convention that
$e_r(x_1,\dots,x_n) = h_r(x_1,\dots,x_n) = 0$ for $r < 0$
and $e_0(x_1,\dots,x_n) = h_0(x_1,\dots,x_n) = 1$.
The following basic identity \cite[$(\mathrm{I}.2.6')$]{Mac} will be used repeatedly:
\begin{equation}\label{mac1}
\sum_{s=0}^r (-1)^{s} e_{s}(x_1,\dots,x_n) h_{r-s}(x_1,\dots,x_n) = 0
\end{equation}
for all $r \geq 1$.
For $m,n \geq 0$, we 
obviously have that
\begin{align}\label{brad1}
h_r(x_1,\dots,x_m,y_1,\dots,y_n) &= \sum_{s=0}^r h_{s}(x_1,\dots,x_m) h_{r-s}(y_{1},\dots,y_n),\\
e_r(x_1,\dots,x_m,y_1,\dots,y_n) &= \sum_{s=0}^r e_{s}(x_1,\dots,x_m) e_{r-s}(y_{1},\dots,y_n).\label{brad12}
\end{align}
Moreover, for all $r \geq 0$, we have that
\begin{align}\label{brad23}
h_r(y_{1},\dots,y_n)&=\sum_{s=0}^r (-1)^{s} e_{s}(x_1,\dots,x_m) 
h_{r-s}(x_{1},\dots,x_m,y_1,\dots,y_n),\\
e_r(y_{1},\dots,y_n)&=\sum_{s=0}^r (-1)^{s} 
h_{s}(x_{1},\dots,x_m)e_{r-s}(x_1,\dots,x_m,y_1,\dots,y_n).\label{brad2}
\end{align}
These identities 
follow by expanding the right hand
sides using the identities (\ref{brad1})--(\ref{brad12}) 
and then simplifying using  (\ref{mac1}).

Let $P$ denote the polynomial algebra
$\C[x_1,\dots,x_n]$, graded so that each $x_i$ is in degree $2$.
The  algebra of {\em symmetric polynomials} is the
subalgebra $P^{S_n}$ of all $S_n$-invariants in $P$.
It is classical that $P$ is a free $P^{S_n}$-module
of rank $n!$
with basis
\begin{equation}\label{bad}
\{x_1^{r_1} x_2^{r_2} \cdots x_n^{r_n}\mid 0 \leq 
r_i < i\}.
\end{equation}
An equivalent statement
is that $P^{S_{n-1}}$
is a free $P^{S_n}$-module
with basis $1,x_n,x_n^2,\dots,x_n^{n-1}$.
The expansion of higher powers of $x_n$ in terms of this basis
can be obtained by
setting $u = x_n$ in the following identity:
for any $r \geq 0$ and any $u$ 
such that
$(u-x_1) \cdots (u-x_n) = 0$ we have that
\begin{equation}\label{mac30}
u^{n+r} = \sum_{s=1}^n  
\sum_{t=0}^r (-1)^{s+t-1} e_{s+t}(x_1,\dots,x_n) h_{r-t}(x_1,\dots,x_n)u^{n-s}.
\end{equation}
This follows by induction on $r$;
the base case $r=0$ is exactly
the assumption 
$(u-x_1) \cdots (u-x_n) = 0$, and then the induction step is obtained by multiplying both sides by $u$ then making some easy manipulations using (\ref{mac1}).
Applying (\ref{mac1}) once more, an equivalent way of writing 
the same identity is as
\begin{equation}\label{mac3}
u^{n+r} = \sum_{s=1}^n\sum_{t=1}^s (-1)^{s-t} e_{s-t}(x_1,\dots,x_n) h_{r+t}(x_1,\dots,x_n) u^{n-s},
\end{equation}
again valid
for all $r \geq 0$ and $u$ 
satisfying
$(u-x_1) \cdots (u-x_n) = 0$

Next fix a composition $\nu$ of $n$ and 
let $P_\nu$ denote the subalgebra $P^{S_\nu}$ of $S_\nu$-invariants in $P$.
For $i \in \Z$
we have the elements 
\begin{align}
h_r(\nu;i) &:= h_r(x_{j+1},x_{j+2},\dots,x_{j+\nu_i}),\\
e_r(\nu;i) &:= e_r(x_{j+1},x_{j+2},\dots,x_{j+\nu_i}),
\end{align}
where $j := \sum_{h < i} \nu_h$.
The algebra $P_\nu$ is freely generated by either 
the $h_r(\nu;i)$'s or the $e_r(\nu;i)$'s
for all $i \in \Z$ and $1 \leq r \leq \nu_i$.
Also, $e_r(\nu;i) = 0$ for $r > \nu_i$.
More generally, given distinct integers $i_1,\dots,i_m$,
let
\begin{align}
h_r(\nu;i_1,\dots,i_m) &:= \sum_{r_1+\cdots+r_m = r}
h_{r_1}(\nu;i_1)h_{r_2}(\nu;i_2) \cdots h_{r_m}(\nu;i_m),\\
e_r(\nu;i_1,\dots,i_m) &:= \sum_{r_1+\cdots+r_m = r}
e_{r_1}(\nu;i_1)e_{r_2}(\nu;i_2) \cdots e_{r_m}(\nu;i_m).
\end{align}
Because of (\ref{brad1})--(\ref{brad12}), 
these are the same elements as defined in the introduction.
If $i_1,\dots,i_m, j_1,\dots,j_l$ are distinct integers
such that $\nu_{i_1}+\cdots+\nu_{i_m}+\nu_{j_1}+\cdots+\nu_{j_l} = n$,
then we have by (\ref{brad23}) that
\begin{align}\label{elo}
h_r(\nu;i_1,\dots,i_m) &= \sum_{s=0}^r (-1)^s
e_s(\nu;j_1,\dots,j_l) h_{r-s}(x_1,\dots,x_n).
\end{align}

The {\em coinvariant algebra} $C$ is the quotient 
$P / I$, where $I$ denotes
the ideal of $P$ generated by all homogeneous symmetric functions of positive degree. 
The images of the monomials (\ref{bad}) give a basis for $C$ as a vector space,
so $\dim C = n!$. 
In fact, 
by a theorem of Chevalley \cite{Chevalley}, the algebra
$C$ viewed as a module over the symmetric group
is isomorphic 
to the left regular module $\C S_n$,
where the action of $S_n$
on $C$ is the action induced by the natural permutation action on $P$.
More generally, 
let $I_\nu$ be the ideal of $P_\nu$ generated by
all homogeneous symmetric polynomials of positive degree
and define the {\em partial coinvariant algebra}
$C_\nu$ to be the quotient $P_\nu / I_\nu$.
The first lemma is well known.

\begin{Lemma} 
The map $P_\nu \rightarrow C$ obtained by restricting the canonical quotient map
$P \twoheadrightarrow C$ to the subalgebra $P_\nu$
has
kernel $I_\nu$ and image $C^{S_\nu}$. Hence it induces a canonical isomorphism between
$C_\nu = P_\nu / I_\nu$ and $C^{S_\nu}$. In particular, 
$\dim C_\nu = |S_n / S_\nu|$.
\end{Lemma}

\begin{proof}
Since we are over a field of characteristic zero, taking $S_\nu$-fixed points is an exact functor.
Applying it to $0 \rightarrow I \rightarrow P \rightarrow C \rightarrow 0$
gives a short exact sequence $0 \rightarrow I^{S_\nu} \rightarrow P_\nu \rightarrow C^{S_\nu}
\rightarrow 0$. So to prove the first statement of the lemma, we just need to show 
that $I^{S_\nu} = I_\nu$.
Let $\gamma: P \twoheadrightarrow P_\nu$ denote the 
projection defined by $\gamma(f) := \frac{1}{|S_\nu|}
\sum_{w \in S_\nu} w f$.
Any element of $I$ is a linear combination of terms of the form
$f g$ for $f \in P$ and $g \in P^{S_n}$ homogeneous of positive degree.
Applying $\gamma$, we  deduce that any element of $\gamma(I)$ is a linear combination of terms of the
form $\gamma(fg) = \gamma(f) g$ for $f \in P$ and $g \in P^{S_n}$
homogeneous of positive degree.
Since $\gamma(f) \in P_\nu$, all such terms belong to $I_\nu$.
So we get that
$I^{S_\nu} = \gamma(I^{S_\nu})
\subseteq \gamma(I) \subseteq I_\nu$.
On the other hand, we obviously have that $I_\nu \subseteq I^{S_\nu}$.
Hence, $I_\nu = I^{S_\nu}$ and
$C_\nu \cong C^{S_\nu}$.
The final statement about dimension now follows easily using Chevalley's theorem, 
since $\dim C^{S_\nu}= \dim (\C S_n)^{S_\nu} = |S_n / S_\nu|$.
\end{proof}

In view of the above lemma,
we will always identify $C_\nu$ with the
subalgebra $C^{S_\nu}$ of $C$.
Moreover, we will use the same notation for elements of $P_\nu$
and for their canonical images in $C_\nu$; since we will usually 
be working in $C_\nu$ from now on this should not cause any confusion.
By (\ref{elo}), we have in $C_\nu$ that
\begin{equation}\label{jack}
h_r(\nu;i_1,\dots,i_m) = (-1)^r e_r(\nu;j_1,\dots,j_l)
\end{equation}
for $m,l \geq 0$ and distinct
integers $i_1,\dots,i_m, j_1,\dots,j_l$
with the property 
that $\nu_{i_1}+\cdots+\nu_{i_m}+\nu_{j_1}+\cdots+\nu_{j_l} = n$.
Hence by (\ref{brad1}) we get that
\begin{equation}\label{jack2}
\sum_{s=0}^r (-1)^s e_s(\nu;i_1,\dots,i_m) h_{r-s}(\nu;i_1,\dots,i_m) = 0
\end{equation}
for all $r \geq 1$, equality again written in $C_\nu$.

Now we come to the {crucial definition}.
Fix another composition $\mu$ of $n$ and let $\lambda$ denote the transpose partition.
Let $I_\nu^\mu$ be the ideal of $P_\nu$ generated by the elements
\begin{equation}\label{gary}
\left\{h_r(\nu;i_1,\dots,i_m)\:\bigg|
\begin{array}{l}
m \geq 1, \:\:i_1,\dots,i_m\text{ distinct integers},\\
r > \lambda_1+\cdots+\lambda_m - \nu_{i_1}-\cdots-\nu_{i_m}\end{array}\right\}
\end{equation}
and set 
\begin{equation}\label{seitz}
C^\mu_\nu := P_\nu / I^\mu_\nu
\end{equation}
exactly as in the introduction. 
If we choose $m$ and $i_1,\dots,i_m$
so that $\lambda_1+\cdots+\lambda_m = \nu_{i_1}+\cdots+\nu_{i_m} = n$,
then $h_r(\nu;i_1,\dots,i_m)$
belongs to
$I^\mu_\nu$ for all $r > 0$,
and it 
equals $h_r(x_1,\dots,x_n)$.
Since $I_\nu$ is generated by the elements
$h_r(x_1,\dots,x_n)$ for all $r > 0$, this shows that
$I_\nu \subseteq I^\mu_\nu$. So it is natural
to regard $C^\mu_\nu$ also as a quotient of $C_\nu$.
Using Lemma~\ref{l2} below, it is easy to see 
that if $\mu$ is regular, i.e. $\lambda_1=n$,
then $I_\nu = I_\nu^\mu$.
So in this special case we have simply that $C^\mu_\nu = C_\nu$.

\begin{Lemma}\label{l2} The ideal
$I^\mu_\nu$ is also generated by
\begin{equation*}
\left\{e_r(\nu;i_1,\dots,i_m)\:\Bigg|
\begin{array}{l}
m \geq 1, \:\: i_1,\dots,i_m\text{ distinct integers},\\
r > \nu_{i_1}+\cdots+\nu_{i_m} - \lambda_{l+1}-\lambda_{l+2}-\cdots\\
\text{where }l:=\#\{i \in \Z\mid \nu_i > 0, i \neq i_1,\dots,i_m\}
\end{array}\right\}.
\end{equation*}
\end{Lemma}

\begin{proof}
Let $J_\nu^\mu$ be the ideal generated by the given set
of elementary symmetric functions.
We just explain how to prove that each generator of $J^\mu_\nu$
belongs to $I^\mu_\nu$, hence $J^\mu_\nu \subseteq I^\mu_\nu$.
Then a similar argument in the other direction
shows that each generator of $I^\mu_\nu$
belongs to $J^\mu_\nu$, hence $I^\mu_\nu \subseteq J^\mu_\nu$,
to complete the proof.

So take some element
$e_r(\nu;i_1,\dots,i_m)$ for $m \geq 1$,
distinct integers $i_1,\dots,i_m$ 
 and
$r > 
\nu_{i_1}+\cdots+\nu_{i_m} - \lambda_{l+1}-\lambda_{l+2}-\cdots$
where $$
l :=\#\{i \in \Z \mid \nu_i > 0, i \neq i_1,\dots,i_m\}.
$$
The definition of $l$
means we can find distinct integers
$j_1,\dots,j_l \notin \{i_1,\dots,i_m\}$ such that 
$\nu_{i_1}+\cdots+\nu_{i_m} + \nu_{j_1}+\cdots+\nu_{j_l} = n$.
It is then the case that
$$
\nu_{i_1}+\cdots+\nu_{i_m} - \lambda_{l+1}-\lambda_{l+2}-\cdots
=
\lambda_1+\cdots+\lambda_l - \nu_{j_1}-\cdots-\nu_{j_l},
$$
hence we have that $r > \lambda_1+\cdots+\lambda_l-\nu_{j_1}-\cdots-\nu_{j_l}$.
If $l > 0$ this shows that
$h_r(\nu;j_1,\dots,j_l)$ is an element of the set (\ref{gary}),
while if $l = 0$ then $r > 0$ so 
$h_r(\nu;j_1,\dots,j_l) = 0$. Either way, this means that
$h_r(\nu;j_1,\dots,j_l) \in I^\mu_\nu$.
We know  by (\ref{jack}) that
$e_r(\nu;i_1,\dots,i_m)$ is equal to
$(-1)^r h_r(\nu;j_1,\dots,j_l)$ plus some element of $I_\nu$.
We observed already above that $I_\nu \subseteq I_\nu^\mu$,
so we have now proved that $e_r(\nu;i_1,\dots,i_m) \in I_\nu^\mu$.
\end{proof}

\begin{Lemma}\label{van}
$C^{\mu}_\nu \neq 0$ if and only if $\lambda \geq \nu^+$.
\end{Lemma}

\begin{proof}
Since everything is graded, we have that $C^\mu_\nu \neq 0$
if and only if all the generators of $I^\mu_\nu$ are of positive
degree, i.e. $\lambda_1+\cdots+\lambda_m \geq \nu_{i_1}+\cdots+\nu_{i_m}$
for all $m \geq 1$ and distinct integers $i_1,\dots,i_m$.
By the definition of the dominance ordering on partitions, this
is the statement that $\lambda \geq \nu^+$.
\end{proof}

When $\nu$ is regular, we write $I^\mu$ for $I^\mu_\nu$ and 
$C^\mu$ for $C^\mu_\nu$.
As we said already in the introduction, Lemma~\ref{l2}
is all that is needed to see that our definition of 
$C^\mu$ is equivalent to
Tanisaki's presentation \cite{T} for the cohomology algebra 
$H^*(F^\mu, \C)$ of the Springer
fiber $F^\mu$. In the following theorem, we record some known 
facts about this algebra.

\begin{Theorem}\label{useful}
Let $\mu$ be a composition of $n$ with transpose partition $\lambda$.
\begin{itemize}
\item[(i)] As a $\C S_n$-module,
$C^\mu$ is isomorphic to the permutation representation $M^\mu$
of $S_n$ on the cosets of the parabolic subgroup $S_\mu$.
\item[(ii)]
The top graded component of $C^\mu$ is in degree
\begin{equation}
d^\mu := \sum_{i \geq 1} \lambda_i(\lambda_i-1)
\end{equation}
(which is twice the dimension of the Springer fiber $F^\mu$).
\item[(iii)] As a $\C S_n$-module,
the top graded component $C^\mu(d^\mu)$ 
is isomorphic to the irreducible Specht module
parametrized by the partition $\mu^+$.
\item[(iv)] 
For any non-zero vector $z \in C^\mu$, there exists
$y \in C^\mu$ such that $yz$ is a non-zero vector
in the top graded component $C^\mu(d^\mu)$.
\end{itemize}
\end{Theorem}

\begin{proof}
Parts (i)--(iii) are proved in \cite{T}.
Part (iv)
is noted in the proof of \cite[Theorem 6.6(vi)]{G}, where it is proved using the monomial basis for $C^\mu$ constructed in
\cite{GP}.
\end{proof}

\begin{Remark}\label{hr}
\rm
Recall from \cite[(I.6.6(vi))]{Mac} that the composition multiplicity
$[M^\mu:S^\tau]$ of the irreducible Specht module $S^\tau$ parametrized by a partition $\tau$ in the permutation module $M^\mu$
is equal to the Kostka number
$K_{\tau,\mu}$.
In view of Theorem~\ref{useful}(i), the polynomial
\begin{equation}
K_{\tau,\mu}(t) := \sum_{r \geq 0} [C^\mu(d^{\mu}-2r): S^\tau]t^r 
\end{equation}
arising from
graded composition multiplicities in $C^\mu$ satisfies 
$K_{\tau, \mu}(1) = K_{\tau,\mu}$.
In fact, it is known by \cite{GP} that $K_{\tau,\mu}(t)$ is
equal to the {\em Kostka-Foulkes polynomial} as defined in 
\cite[section III.6]{Mac}.
\end{Remark}

\section{Partial coinvariant algebras}\label{pca}

In the remainder of the article, the notation $\bigoplus_\nu$ always denotes
the direct sum over the set of all compositions $\nu$ of $n$.
Let $\hat{\mathfrak{g}} := \mathfrak{gl}_\infty(\C)$
be the Lie algebra 
of matrices over $\C$ with rows and columns parametrized by index set $\Z$,
all but finitely many of whose entries are
zero. For each $i \in \Z$, let ${D}_i, {E}_i$ and 
$F_i$ denote the $(i,i)$-matrix unit, the $(i,i\!+\!1)$-matrix unit 
and the $(i\!+\!1,i)$-matrix unit 
in $\hat{\mathfrak{g}}$, respectively.
These elements 
generate $\hat{\mathfrak{g}}$.

We are going to exploit some basic facts about
{\em polynomial representations} of $\hat{\mathfrak{g}}$.
All of these facts are standard results about polynomial representations
finite dimensional general linear Lie algebras
(see e.g. \cite{Green})
extended to $\hat{\mathfrak{g}}$ by 
taking direct limits.
To start with, given a $\hat{\mathfrak{g}}$-module $M$ and 
a composition $\nu$ of $n$, the {\em $\nu$-weight space}
of $M$ is
\begin{equation}
M_\nu := \{v \in M\mid D_i v = \nu_i v \text{ for all }i\in\Z\}.
\end{equation}
We call $M$ a {\em polynomial representation} of $\hat{\mathfrak{g}}$
of {\em degree} $n$ if 
$M= \bigoplus_\nu M_\nu$
and all $M_\nu$ are finite dimensional.
For example, the $n$th tensor power $\widehat{V}^{\otimes n}$
of the natural $\hat{\mathfrak{g}}$-module $\widehat{V}$
is a polynomial representation of
degree $n$, as is the module
\begin{equation}
{\textstyle\bigwedge}^\mu(\widehat{V}) := \cdots \otimes {\textstyle\bigwedge}^{\mu_1}(\widehat{V})
\otimes {\textstyle\bigwedge}^{\mu_2}(\widehat{V}) \otimes\cdots
\end{equation}
for a composition $\mu$ of $n$.
For each partition $\lambda$ of $n$, there is an
irreducible polynomial representation $P^\lambda(\widehat{V})$
characterized uniquely up to isomorphism by the property that
the dimension of the $\nu$-weight space 
of $P^\lambda(\widehat{V})$ is equal to the Kostka number $K_{\lambda,\nu}$
 for every $\nu$.
These modules give all of the irreducible polynomial representations of degree $n$.

The symmetric group $S_n$ acts on the right on 
$\widehat{V}^{\otimes n}$ by permuting tensors.
Consider the functor 
$\widehat{V}^{\otimes n} \otimes_{\C S_n} ?$ from the category of finite 
dimensional left
$\C S_n$-modules to the category of polynomial representations of
$\widehat{\mathfrak{g}}$ of degree $n$. By Schur's classical theory,
this functor is known to be an equivalence of categories.
Recall that $M^\mu$ is the permutation module parametrized by a composition $\mu$ 
of $n$ 
and $S^\lambda$ is the irreducible
Specht module parametrized by a partition $\lambda$ of $n$.
It is well known that
\begin{align}\label{hand1}
{\textstyle\bigwedge}^{\mu}(\widehat{V}) &\cong \widehat{V}^{\otimes n} \otimes_{\C S_n} \widetilde{M}^{\mu},\\
P^{\lambda}(\widehat{V}) &\cong \widehat{V}^{\otimes n} \otimes_{\C S_n} 
S^{\lambda},\label{hand2}
\end{align}
where 
$\widetilde{M}^\mu$ denotes the $\C S_n$-module
obtained from $M^\mu$ by twisting the action by sign.

\begin{Lemma}\label{b4}
For any left $\C S_n$-module $M$,
there is a natural vector space isomorphism
$\kappa_M:\widehat{V}^{\otimes n} \otimes_{\C S_n} M
\stackrel{\sim}{\rightarrow} \bigoplus_\nu M^{S_\nu}$ defined as follows.
Suppose $\nu$ is a composition of $n$ and 
$a_1 \leq \cdots \leq a_n$ are integers
exactly $\nu_i$ of which
are equal to $i$ for each $i \in \Z$.
Then 
\begin{align*}
\kappa_M(v_{a_1} \otimes \cdots \otimes v_{a_n} \otimes m)&= \frac{1}{|S_\nu|} \sum_{w \in S_\nu} w m \in M^{S_\nu}\\\intertext{for any $m \in M$.
The inverse map $\kappa_M^{-1}$ satisfies}
\kappa_M^{-1}(m) &= 
v_{a_1}\otimes\cdots\otimes v_{a_n} \otimes m
\end{align*}
for any $m \in M^{S_\nu}$.
\end{Lemma}

\begin{proof}
Fix a composition $\nu$ of $n$
and let $a_1 \leq \cdots \leq a_n$ be the unique integers
such that exactly $\nu_i$ of them are equal to $i$ for each $i \in \Z$.
Let $\widehat{V}^{\otimes n}_\nu$ denote the $\nu$-weight space
of $\widehat{V}^{\otimes n}$.
Writing $\C$ for the trivial module, 
it is well known that
the right $\C S_\nu$-module homomorphism 
$\C \rightarrow \widehat{V}^{\otimes n}_\nu$
under which
$1 \mapsto v_{a_1}\otimes\cdots\otimes v_{a_n}$
extends by Frobenius reciprocity to a right $\C S_n$-module
isomorphism 
$\C \otimes_{\C S_\nu} \C S_n \stackrel{\sim}{\rightarrow}
\widehat{V}^{\otimes n}_\nu$. 
Also there is a familiar vector space isomorphism 
$\C \otimes_{\C S_\nu} M \stackrel{\sim}{\rightarrow} M^{S_\nu},
1 \otimes m \mapsto \frac{1}{|S_\nu|} \sum_{w \in S_\nu} w m$.
Composing these maps, we obtain an isomorphism
$$
\widehat{V}^{\otimes n}_\nu \otimes_{\C S_n} M
\stackrel{\sim}{\longrightarrow}
\C \otimes_{\C S_\nu} \C S_n \otimes_{\C S_n} M
= 
\C \otimes_{\C S_\nu} M
\stackrel{\sim}{\longrightarrow} M^{S_\nu}
$$
such that $v_{a_1}\otimes\cdots\otimes v_{a_n} \otimes m
\mapsto \frac{1}{|S_\nu|} \sum_{w \in S_\nu} wm$.
The isomorphism $\kappa_M$
is the direct sum of these maps
over all compositions $\nu$ of $n$.
\end{proof}

As in the introduction, 
let $F$ be the flag manifold of (complex)
dimension $\frac{1}{2}n(n-1)$, and
let $F_\nu$ denote the partial flag manifold associated to a composition $\nu$ of $n$.
So elements of $F_\nu$ are nested 
chains $(V_j)_{j \in \Z}$ of subspaces of $\C^n$ such that
$\dim V_j = \sum_{i \leq j} \nu_i$ for each $j$.
The (complex) dimension of $F_\nu$ is equal to $\frac{1}{2}d_\nu$ where
\begin{equation}\label{dnu}
d_\nu := n(n-1) - \sum_{i \in\Z} \nu_i(\nu_i-1).
\end{equation}
Let $\pi:F \twoheadrightarrow F_\nu$ denote the natural projection;
informally, $\pi$ is the map forgetting 
all subspaces of a full flag of the wrong dimension.
Identifying the cohomology algebra $H^*(F, \C)$ 
with the coinvariant algebra $C$ 
as in \cite[10.2(3)]{Fulton}, the fundamental class of a point of $F$ 
(regarded as an element of $H^*(F,\C)$ via Poincar\'e duality)
is the canonical image of the polynomial
\begin{equation}
\eps := \frac{1}{n!}\prod_{1 \leq i < j \leq n} (x_i - x_j)
\end{equation}
in $C$.
More generally, the 
cohomology algebra $H^*(F_\nu, \C)$ can be identified 
with the partial coinvariant algebra $C_\nu = C^{S_\nu}$
so that the pull-back homomorphism
$\pi^*: H^*(F_\nu, \C) \rightarrow H^*(F, \C)$
coincides with the
natural inclusion $C_\nu \hookrightarrow C$.
Let
\begin{equation}\label{epsnu}
\eps_\nu := \frac{1}{|S_\nu|}
\prod_{\substack{1 \leq i < j \leq n \\
i\: \stackrel{\nu}{\sim}\: j}} (x_i - x_j),
\end{equation}
where we write 
$i \stackrel{\nu}{\sim} j$ if $i$ and $j$ lie in the same $S_\nu$-orbit.
Note $\eps$ is divisible by $\eps_\nu$ 
and the quotient 
$\eps / \eps_\nu$ belongs to $P_\nu$.
The fundamental class of a point of $F_\nu$
is the canonical image of $\eps / \eps_\nu$ in $C_\nu$.

Let $\widetilde{C}$ denote the $\C S_n$-module obtained from $C$
by twisting the natural permutation action by sign, i.e. $w \in S_n$ acts on $\widetilde{C}$ as the map
$f \mapsto \sgn(w) wf$.
Since the regular $\C S_n$-module 
twisted by sign is 
isomorphic to the regular module, it follows from Chevalley's theorem that
$\widetilde{C}$ is isomorphic to $C$ as (ungraded) 
$\C S_n$-modules.
Let $\widetilde{C}_\nu$ denote the space $\widetilde C^{S_\nu}$
of all $S_\nu$-invariants in $\widetilde{C}$. In other words, 
$\widetilde{C}_\nu$ is the space of all $S_\nu$-{\em anti-invariants} in $C$. 
Note $\eps_\nu$ belongs to
$\widetilde{C}_\nu$.

In the next lemma, we use Poincar\'e duality to 
regard the usual push-forward 
$\pi_*$ in homology as a homogeneous map
$\pi_*:H^*(F,\C) \rightarrow H^*(F_\nu,\C)$
of degree $-\sum_{i \in \Z} \nu_i(\nu_i-1)$.

\begin{Lemma}\label{b3}
The restriction of the push-forward 
$\pi_*:H^*(F, \C) \rightarrow H^*(F_\nu, \C)$ 
to the subspace $\widetilde{C}_\nu$
defines a $C_\nu$-module isomorphism
$\pi_*:\widetilde{C}_\nu\stackrel{\sim}{\rightarrow}C_\nu$
with $\pi_*(\eps_\nu) = 1$.
In particular, $\widetilde{C}_\nu = \eps_\nu C_\nu$
and $\dim \widetilde{C}_\nu = |S_n / S_\nu|$.
\end{Lemma}

\begin{proof}
By the projection formula, $\pi_*$ is a $C_\nu$-module homomorphism, 
hence so is its restriction to the $C_\nu$-submodule
$\widetilde{C}_\nu$ of $C$.
By degree considerations,
we have that $\pi_*(\eps_\nu) = c\cdot 1$
for some scalar $c$.
To compute the scalar,
$\pi_*$ sends the fundamental class of a point to the fundamental class of a point,
so $\pi_*(\eps) = \eps / \eps_\nu$.
Hence in $C_\nu$ we have that
$$
\eps / \eps_\nu = \pi_*(\eps) =
\pi_*(\eps_\nu \cdot \eps / \eps_\nu)
= \pi_*(\eps_\nu) \eps / \eps_\nu = c \eps / \eps_\nu.
$$
So $c=1$ and $\pi_*(\eps_\nu) = 1$.
It follows at once from this that
$\pi_*: \widetilde{C}_\nu \rightarrow C_\nu$ is surjective.
It is an isomorphism because $\dim \widetilde{C}_\nu = \dim C_\nu = |S_n / S_\nu|$.
\end{proof}

Now let us consider the $\hat{\mathfrak{g}}$-module
$\widehat{V}^{\otimes n} \otimes_{\C S_n} \widetilde{C}$.
Since $\widetilde{C} \cong \C S_n$ as a $\C S_n$-module,
this is isomorphic simply to tensor space $\widehat{V}^{\otimes n}$, but it carries an interesting grading induced from the natural grading on $C$.
Composing
the direct sum 
over all $\nu$ of 
the isomorphisms from Lemma~\ref{b3}
with
the isomorphism $\kappa_{\widetilde{C}}$ from Lemma~\ref{b4}, 
we obtain an isomorphism
\begin{equation}\label{phi}
\phi: \widehat{V}^{\otimes n} \otimes_{\C S_n} \widetilde C
\stackrel{\sim}{\longrightarrow} {\textstyle\bigoplus_\nu} C_\nu.
\end{equation}
Using this, we transport the natural 
$\hat{\mathfrak{g}}$-action
on $\widehat{V}^{\otimes n} \otimes_{\C S_n} \widetilde C$
to the vector space $\bigoplus_\nu C_\nu$, to make the latter into 
a $\hat{\mathfrak{g}}$-module whose 
$\nu$-weight space is equal to the partial coinvariant algebra 
$C_\nu$.
Because this definition involves push-forward,
the action of $\hat{\mathfrak{g}}$ 
does not preserve the grading on $\bigoplus_\nu C_\nu$ in the usual sense.
Instead, $\hat{\mathfrak{g}}$ leaves the subspaces
$\bigoplus_\nu C_\nu(d_\nu - 2r)$
invariant for each $r \geq 0$.

Consider the following {\em key situation}.
Fix $i \in \Z$, let $\nu$ be a composition of $n$
with $\nu_i \neq 0$, and define
$\nu'$ to be
the composition of $n$ obtained from $\nu$ by replacing
$\nu_i$ by $\nu_i-1$ and $\nu_{i+1}$ by $\nu_{i+1}+1$.
The Chevalley generators $E_i$ and $F_i$ of $\hat{\mathfrak{g}}$
define linear maps
\begin{equation}\label{maps}
C_\nu
\quad\qquad
C_{\nu'}.
\begin{picture}(0,20)
\put(-35,4){\makebox(0,0){$\stackrel{\stackrel{\,\scriptstyle{F_i}_{\phantom{,}}}{\displaystyle\longrightarrow}}{\stackrel{\displaystyle\longleftarrow}{\scriptstyle{E_i}}}$}}
\end{picture}
\end{equation}
\vspace{0mm}

\noindent
We are going to calculate these maps explicitly.
It will be convenient to let $a := \nu_i'$, 
$b := \nu_{i+1}$ and $k := \sum_{j \leq i} \nu_j$, so that
\begin{align}
&x_{k-a} \stackrel{\nu}{\sim} 
\cdots
\stackrel{\nu}{\sim} x_{k-1}
\stackrel{\nu}{\sim} x_{k},\quad x_{k+1}
\stackrel{\nu}{\sim}
\cdots
\stackrel{\nu}{\sim} x_{k+b},\\
&x_{k-a} \stackrel{\nu'}{\sim} 
\cdots
\stackrel{\nu'}{\sim} x_{k-1},\quad
x_{k}
\stackrel{\nu'}{\sim}  x_{k+1}
\stackrel{\nu'}{\sim}
\cdots
\stackrel{\nu'}{\sim} x_{k+b}.
\end{align}
The notation $\stackrel{\nu}{\sim}$ 
and $\stackrel{\nu'}{\sim}$ being used here 
was introduced after (\ref{epsnu}).

Let
\begin{equation}
F_{\nu,\nu'} :=
\big\{((V_j)_{j \in \Z}, (V_j')_{j \in \Z})
\in F_\nu \times F_{\nu'}\:\big|\:
V_j = V_j'\text{ for }j \neq i,\:
V_i' \subset V_i\big\}.
\end{equation}
This is obviously isomorphic to the partial flag manifold
whose cohomology algebra has been identified with
$C_{\nu,\nu'} := C^{S_\nu \cap S_{\nu'}}$.
Both $C_\nu = C^{S_\nu}$ and $C_{\nu'} = C^{S_{\nu'}}$ are subalgebras
of $C_{\nu,\nu'}$.
By the sentence after (\ref{bad}),
$C_{\nu,\nu'}$ is a free $C_\nu$-module with basis
$1,x_k,\dots,x_k^a$ and a free $C_{\nu'}$-module with basis
$1,x_k,\dots,x_k^b$.
Moreover, by (\ref{mac3}), we have that for $r \geq 1$ that
\begin{align}\label{mac4}
x_k^{a+r} &= \sum_{s=0}^{a}\sum_{t=0}^{s}  (-1)^{s-t}
e_{s-t}(\nu;i) h_{r+t}(\nu;i) x_k^{a-s},\\
\label{mac5}
x_k^{b+r} &= \sum_{s=0}^{b}\sum_{t=0}^{s} 
 (-1)^{s-t} e_{s-t}(\nu';i+1) h_{r+t}(\nu';i+1)
x_k^{b-s}.
\end{align}

Let
$p:F_{\nu,\nu'} \rightarrow F_\nu$ and
$p':F_{\nu,\nu'} \rightarrow F_{\nu'}$ be the first and second
projections.
The pull-backs $p^*:C_\nu \rightarrow C_{\nu,\nu'}$
and $(p')^*:C_{\nu'} \rightarrow C_{\nu,\nu'}$ are simply the natural 
inclusions.
To describe the push-forwards (again regarded as maps in cohomology
via Poincar\'e duality),
we define
\begin{equation}\label{epsnunu}
\eps_{\nu,\nu'} := \frac{1}{|S_\nu \cap S_{\nu'}|}
\prod_{\substack{1 \leq i < j \leq n \\
i\: \stackrel{\nu}{\sim}\: j
\: \stackrel{\,\nu{\scriptscriptstyle'}}{\!\sim}\: i}} (x_i - x_j)
\end{equation}
like in (\ref{epsnu}),
so the fundamental class of a point of $F_{\nu,\nu'}$
is the canonical image of 
$\eps / \eps_{\nu,\nu'}$ in $C_{\nu,\nu'}$.
Both $\eps_\nu$ and $\eps_{\nu'}$ are divisible by $\eps_{\nu,\nu'}$, and
we have that
\begin{align}\label{boris1}
\eps_\nu / \eps_{\nu,\nu'} &= 
\frac{1}{a\!+\!1} \prod_{j=1}^a (x_{k-j}-x_k)
=
\frac{(-1)^a}{a\!+\!1}
\sum_{s=0}^a (-1)^s e_s(\nu';i) x_k^{a-s},
\\
\eps_{\nu'} / \eps_{\nu,\nu'} &= 
\frac{1}{b\!+\!1} \prod_{j=1}^b (x_k-x_{k+j})
=
\frac{1}{b\!+\!1}
\sum_{s=0}^b (-1)^s e_s(\nu;i+1) x_k^{b-s}.\label{boris2}
\end{align}

\begin{Lemma}\label{push}
Let notation be as above, so in particular $a =\nu_i-1 = \nu_{i}'$ and 
$b = \nu_{i+1} = \nu_{i+1}' - 1$.
\begin{itemize}
\item[(i)]
The push-forward $p_*:C_{\nu,\nu'} \rightarrow C_\nu$
is the unique homogeneous
$C_\nu$-module homomorphism 
of degree 
$-a$ that maps $\eps_\nu / \eps_{\nu,\nu'}$ to $1$.
Equivalently, it is the unique 
$C_\nu$-module homomorphism with
$p_*(x_k^r) = (-1)^a h_{r-a}(\nu;i)$ for each $r \geq 0$.
\item[(ii)]
The push-forward $p'_*:C_{\nu,\nu'} \rightarrow C_{\nu'}$
is the unique homogeneous 
$C_{\nu'}$-module homomorphism of degree $-b$
that maps
$\eps_{\nu'} / \eps_{\nu,\nu'}$ to $1$.
Equivalently, it is the unique $C_{\nu'}$-module homomorphism 
with $p'_*(x_k^r)=h_{r-b}(\nu';i+1)$ for each $r \geq 0$.
\end{itemize}
\end{Lemma}

\begin{proof}
By the projection formula, $p_*$ is a homogeneous $C_\nu$-module
homomorphism of degree $-a$,
so it must map $\eps_\nu / \eps_{\nu,\nu'}$ to $c\cdot 1$ for some
scalar $c$.
Since it maps the fundamental class of a point to the fundamental
class of a point, we know that
$p_*(\eps / \eps_{\nu,\nu'}) = \eps / \eps_\nu$, from which we get
that $c = 1$ as in the proof of Lemma~\ref{b3}.
Now suppose that $f:C_{\nu,\nu'} \rightarrow C_\nu$
is any homogeneous $C_\nu$-module homomorphism of degree $-a$
mapping $\eps_\nu / \eps_{\nu,\nu'}$ to $1$.
By (\ref{boris1}) and (\ref{brad2}), 
\begin{align*}
\eps_\nu / \eps_{\nu,\nu'} &=
\frac{(-1)^a}{a\!+\!1}
\sum_{s=0}^a\sum_{t=0}^s (-1)^{s+t} e_{s-t}(\nu;i)
 x_k^{a-s+t}.
\end{align*}
Applying $f$ to this equation 
and observing that
$f(1) = \cdots =f(x_k^{a-1}) = 0$ by degree considerations,
we deduce that $f(x_k^a) = (-1)^a$. Finally using (\ref{mac4}),
we get that $f(x_k^r) = (-1)^a h_{r-a}(\nu;i)$ for any 
$r \geq 0$ (recalling 
by convention that $h_{r-a}(\nu;i) = 0$ if $r < a$).
Since the elements $x_k^r$ generate $C_{\nu,\nu'}$ as a $C_\nu$-module, there is clearly a unique such map.
This proves (i), and (ii) is similar.
\end{proof}

\begin{Theorem}\label{be}
Let notation be as above, so in particular $a =\nu_i-1 = \nu_{i}'$ and 
$b = \nu_{i+1} = \nu_{i+1}' - 1$.
\begin{itemize}
\item[(i)]
The map $F_i:C_\nu \rightarrow C_{\nu'}$ is equal to
the composite
$$
C_\nu \stackrel{p^*}{\longrightarrow}
C_{\nu,\nu'}\stackrel{m}{\longrightarrow}
C_{\nu,\nu'} \stackrel{p'_*}{\longrightarrow} C_{\nu'}
$$
where $m$ is the map defined by multiplication by $\prod_{j=1}^a (x_{k-j}-x_k)$.
Equivalently, it is the restriction to $C_\nu$ of the unique $C_{\nu'}$-module homomorphism $C_{\nu,\nu'} \rightarrow C_{\nu'}$ sending
$$
x_k^r \mapsto (-1)^a\sum_{s=0}^a (-1)^s e_{s}(\nu';i) h_{r-s+a-b}(\nu';i+1)
$$
for each $r \geq 0$.

\item[(ii)]
The map $E_i:C_{\nu'} \rightarrow C_{\nu}$ is equal to
the composite
$$
C_{\nu'} \stackrel{(p')^*}{\longrightarrow} C_{\nu,\nu'}
\stackrel{m'}{\longrightarrow} C_{\nu,\nu'}
\stackrel{p_*}{\longrightarrow} C_{\nu}
$$
where $m'$ is the map defined by 
multiplication by $\prod_{j=1}^b (x_k-x_{k+j})$.
Equivalently, it is the restriction to $C_{\nu'}$ 
of the unique $C_{\nu}$-module homomorphism $C_{\nu,\nu'} \rightarrow C_{\nu}$ sending
$$
x_k^r \mapsto 
(-1)^a \sum_{s=0}^b (-1)^{s} e_s(\nu;i+1) h_{r-s+b-a}(\nu;i)
$$
for each $r \geq 0$.
\end{itemize}
\end{Theorem}

\begin{proof}
We just prove (i), the proof of (ii) being similar.
Observe to start with that the action of
$\hat{\mathfrak{g}}$ on $\bigoplus_\nu C_\nu$ lifts 
to an action on $\bigoplus_\nu P_\nu$.
To see this, let $\widetilde{P}$ denote the $\C S_n$-module obtained from $P$
by twisting the natural permutation action by sign.
Let $\widetilde{P}_\nu$ denote the 
$S_\nu$-invariants in $\widetilde{P}$,
or equivalently, the $S_\nu$-anti-invariants in $P$.
By another fundamental theorem of Chevalley, 
every element of $\widetilde{P}_\nu$ is divisible by $\eps_\nu$.
Division by $\eps_\nu$ defines a $P_\nu$-module isomorphism
$\hat\pi_*:\widetilde{P}_\nu \stackrel{\sim}{\rightarrow} P_\nu$
lifting the isomorphism
$\pi_*:\widetilde{C}_\nu \stackrel{\sim}{\rightarrow} C_\nu$
from Lemma~\ref{b3}.
Composing the direct sum of these isomorphisms over all
$\nu$ with the 
isomorphism
$\kappa_{\widetilde{P}}$
from Lemma~\ref{b4}, we obtain an isomorphism
$\hat\phi$
making the following diagram commute:
$$
\begin{CD}
\widehat{V}^{\otimes n} \otimes_{\C S_n} \widetilde{P}
&@>\hat\phi>>& \bigoplus_\nu P_\nu \\
@VVV&&@VVV\\
\widehat{V}^{\otimes n} \otimes_{\C S_n} \widetilde{C}
&@>>\phi >& \bigoplus_\nu C_\nu \\
\end{CD}
$$
where the vertical maps here (and later on)
are the canonical quotient maps.
The natural action of $\hat{\mathfrak{g}}$
on 
$\widehat{V}^{\otimes n} \otimes_{\C S_n} \widetilde P$
clearly lifts the action of $\hat{\mathfrak{g}}$ on 
$\widehat{V}^{\otimes n} \otimes_{\C S_n} \widetilde C$.
Transporting this action through the isomorphism $\hat\phi$
gives the desired action of $\hat{\mathfrak{g}}$ on $\bigoplus_\nu P_\nu$.

Now we compute the effect of
$F_i$ on a polynomial $f \in P_\nu$.
The inverse image of $f$ under $\hat\phi$ is
$(\cdots \otimes v_i^{\otimes (a+1)} \otimes v_{i+1}^{\otimes b}
\otimes\cdots) \otimes \eps_\nu f$.
Applying the Lie algebra element $F_i$ (which maps one $v_i$
to $v_{i+1}$ in all possible ways),
we get
$$
(\cdots \otimes v_i^{\otimes a} \otimes v_{i+1}^{\otimes (b+1)}
\otimes\cdots)(1+\sum_{j=1}^a (k\!-\!j\:\:k))
\otimes \eps_\nu f.
$$
Since $-(k\!-\!j\:\:k) \eps_\nu f = \eps_\nu f$, this
equals 
$(\cdots \otimes v_i^{\otimes a} \otimes v_{i+1}^{\otimes (b+1)}\otimes \cdots)
\otimes (a+1)\eps_\nu f$,
which is the same as
$(\cdots \otimes v_i^{\otimes a} \otimes v_{i+1}^{\otimes (b+1)}\otimes \cdots)
\otimes \eps_{\nu,\nu'} \prod_{j=1}^a (x_{k-j}-x_k) f$
by (\ref{boris1}).
Applying $\hat\phi$, we get
$$
\frac{1}{\eps_{\nu'}}
\frac{1}{|S_{\nu'}|} 
\sum_{w \in S_{\nu'}} \sgn (w) w\left (\eps_{\nu,\nu'} {\textstyle
\prod_{j=1}^a (x_{k-j}-x_k)} f\right).
$$
Setting $P_{\nu,\nu'} := P^{S_\nu \cap S_{\nu'}}$,
we have proved that the map
$F_i:P_\nu \rightarrow P_{\nu'}$ is equal to the composite
of the inclusion $\hat p^*:P_\nu \hookrightarrow P_{\nu,\nu'}$, then
the map $\hat m:P_{\nu,\nu'} 
\rightarrow P_{\nu,\nu'}$
defined by multiplication by 
$\prod_{j=1}^a (x_{k-j}-x_k)$, and finally
the map $\hat p'_*:P_{\nu,\nu'} \rightarrow P_{\nu'}$
defined by
$$
f \mapsto 
\frac{1}{\eps_{\nu'}}
\frac{1}{|S_{\nu'}|}
\sum_{w \in S_{\nu'}} \sgn (w) w\left (\eps_{\nu,\nu'} f\right).
$$
To complete the proof of the first statement of (i), 
it just remains to descend back down to partial coinvariant algebras,
which 
amounts to checking that 
all three squares in the following diagram commute:
$$
\begin{CD}
P_\nu &@>\hat p^*>>& P_{\nu,\nu'}&@>\hat m>>&P_{\nu,\nu'}&@>\hat p'_*>>&P_{\nu'}\phantom{.}\\
@VVV&&@VVV&&@VVV&&@VVV\\
C_\nu &@>>p^*>& C_{\nu,\nu'}&@>>m>&C_{\nu,\nu'}&@>>p'_*>&C_{\nu'}.
\end{CD}
$$
The commutativity of the left hand two squares is obvious.
For the right hand square, 
note from its definition that $\hat p'_*$ is a
homogeneous $P_{\nu'}$-module homomorphism of degree $-b$
mapping $\eps_{\nu'} / \eps_{\nu,\nu'}$ to $1$.
Moreover,
for any $f \in P_{\nu,\nu'}$ and any 
homogeneous 
symmetric function $g$ of positive degree,
we have that $\hat p'_*(fg) = \hat p'_*(f)g$, so $\hat p'_*$ factors
through the quotients to induce a 
homogeneous $C_{\nu'}$-module homomorphism
$C_{\nu,\nu'} \rightarrow C_{\nu'}$ of degree $-b$ mapping $\eps_{\nu'} / \eps_{\nu,\nu'}$ to $1$.
This
coincides with the map $p'_*$ by Lemma~\ref{push}(ii).

Finally, to deduce  the second description of $F_i$,
take an element $z \in C_\nu$ and write 
$z = \sum_{r=0}^b z_r x_k^r$ for (unique) elements $z_r \in C_{\nu'}$.
Using (\ref{boris1}) 
and the second description of $p'_*$ from Lemma~\ref{push}(ii),
$p'_*\circ m \circ p^*$ maps $z$ to
$$
\sum_{r=0}^b z_r \cdot (-1)^a \sum_{s=0}^a (-1)^s e_s(\nu';i)
h_{r-s+a-b}(\nu';i+1).
$$
This is also the image of $z$ under the given $C_{\nu'}$-module
homomorphism $C_{\nu,\nu'}\rightarrow C_{\nu'}$.
\end{proof}

\section{The algebras $C^\mu_\nu$}

Throughout the section, we fix a composition $\mu$ of $n$
with transpose partition $\lambda$.
We will now regard $I^\mu_\nu$ as an ideal
of $C_\nu$ rather than of $P_\nu$, 
generated by the canonical images of the elements (\ref{gary}).
So, recalling (\ref{seitz}), we are now viewing 
the algebra $C^\mu_\nu$ as the quotient $C_\nu / I^\mu_\nu$, and
$\bigoplus_\nu C^\mu_\nu$ is the
quotient of $\bigoplus_\nu C_\nu$ by the subspace $\bigoplus_\nu I^\mu_\nu$.
The following lemma implies that the action of $\hat{\mathfrak{g}}$
on $\bigoplus_\nu C_\nu$ factors through this quotient to induce a well-defined
action of $\hat{\mathfrak{g}}$ on $\bigoplus_\nu C^\mu_\nu$.

\begin{Lemma}\label{dlem}
For each $i \in \Z$, 
the Chevalley generators $E_i$ and $F_i$ of $\hat{\mathfrak{g}}$ leave the subspace
$\bigoplus_\nu I^\mu_\nu$
of $\bigoplus_\nu C_\nu$ invariant.
\end{Lemma}

\begin{proof}
Fix $i \in \Z$ and let $\nu,\nu'$ be compositions as in 
the key situation (\ref{maps}).
We also let $a := \nu_i', b := \nu_{i+1}$ and $k := \sum_{j \leq i} \nu_j$. 
It suffices
to show that
$F_i(I^\mu_\nu) \subseteq I^\mu_{\nu'}$ and
$E_i(I^{\mu}_{\nu'}) \subseteq I^\mu_\nu$.
We just verify the first containment, the second being entirely similar.
Since $C_\nu \subseteq C_{\nu,\nu'}$ and $C_{\nu,\nu'}$ is generated
as a $C_{\nu'}$-module by the elements $x_k^p$ for $p \geq 0$,
every element of
$I_\nu^\mu$ is a $C_{\nu'}$-linear combination of
terms of the form
$x_k^p h_r(\nu;i_1,\dots,i_m)$
for $p \geq 0$, $m \geq 1$, distinct integers $i_1,\dots,i_m$
and $r > \lambda_1+\cdots+\lambda_m - \nu_{i_1}-\cdots-\nu_{i_m}$.
Recalling the second description of $F_i$ from Theorem~\ref{be}(i),
it therefore suffices
to show that $F_i$ maps all these 
$x_k^p h_r(\nu;i_1,\dots,i_m)$'s
into $I^\mu_{\nu'}$.

If neither $i$ nor $(i+1)$ 
belongs to the set $\{i_1,\dots,i_m\}$,
then $h_r(\nu;i_1,\dots,i_m)$ equals $h_r(\nu';i_1,\dots,i_m)\in C_{\nu'}$. So applying Theorem~\ref{be}(i), we see that $F_i$ maps
$x_k^p h_r(\nu;i_1,\dots,i_m)$ to 
$$
(-1)^a \sum_{s=0}^a (-1)^s e_s(\nu';i)
h_{p-s+a-b}(\nu';i+1)
h_r(\nu';i_1,\dots,i_m),
$$
which belongs to $I_{\nu'}^\mu$ because $h_r(\nu';i_1,\dots,i_m)$ does 
already according to the definition (\ref{gary}).
The same argument applies if both $i$ and $(i+1)$ belong to the
set $\{i_1,\dots,i_m\}$.
So it remains to consider the 
situations when exactly one of $i$ and $(i+1)$
belongs to $\{i_1,\dots,i_m\}$.
We may as well assume either that $i_m = i$ or that $i_m = i+1$.

Suppose first that $i_m = i$ and 
$i_1,\dots,i_{m-1} \neq i+1$.
By (\ref{brad1}), we have that
$x_k^p h_r(\nu;i_1,\dots,i_{m}) = \sum_{t=0}^r 
h_t(\nu';i_1,\dots,i_{m}) x_k^{p+r-t}$.
Applying $F_i$ using Theorem~\ref{be}(i), we get
the element
$$
(-1)^a\sum_{s=0}^a\sum_{t=0}^r  (-1)^s e_s(\nu';i)
h_{p+r-t-s+a-b}(\nu';i+1)
h_t(\nu';i_1,\dots,i_{m}).
$$
Since
$r > \lambda_1+\cdots+\lambda_m
- \nu_{i_1}-\cdots-\nu_{i_m} = 
\lambda_1+\cdots+\lambda_m-\nu_{i_1}'-\cdots-\nu_{i_m}'-1$, we have that
$h_t(\nu';i_1,\dots,i_{m}) \in I^\mu_{\nu'}$ if $t > r$
by (\ref{gary}),
so another application of (\ref{brad1}) gives that
\begin{multline*}
\sum_{t=0}^r
h_{p+r-t-s+a-b}(\nu';i+1)
h_t(\nu';i_1,\dots,i_m)\\
\equiv
h_{p+r-s+a-b}(\nu;i_1,\dots,i_m,i+1)
\pmod{I_{\nu'}^\mu}.
\end{multline*}
It remains to show that
$\sum_{s=0}^a (-1)^s e_s(\nu';i)
h_{p+r-s+a-b}(\nu';i_1,\dots,i_{m},i+1)$
belongs to $I_{\nu'}^\mu$.
By (\ref{brad23}), it equals
$h_{p+r+a-b}(\nu';i_1,\dots,i_{m-1},i+1)$
which does indeed lie in $I_{\nu'}^\mu$ since
$p+r+a-b
= p + r + \nu_i-\nu_{i+1}'
> \lambda_1+\cdots+\lambda_m - \nu_{i_1}'-\cdots-\nu_{i_{m-1}}' -\nu_{i+1}'$.

Finally suppose that $i_m = i+1$ and $i_1,\dots,i_{m-1} \neq i$.
By (\ref{brad23}), we can rewrite
$x_k^p h_r(\nu;i_1,\dots,i_{m})$
as
$x_k^p h_r(\nu';i_1,\dots,i_{m})
- x_k^{p+1} h_{r-1}(\nu';i_1,\dots,i_{m})$.
The image of this under $F_i$ is
\begin{multline*}
(-1)^a\sum_{s=0}^a
(-1)^s e_s(\nu';i)
\big(
h_{p-s+a-b}(\nu';i+1) h_r(\nu';i_1,\dots,i_{m})
\\-
h_{p+1-s+a-b}(\nu';i+1) h_{r-1}(\nu';i_1,\dots,i_{m})
\big).
\end{multline*}
Since $r > 
 \lambda_1+\cdots+\lambda_m
- \nu_{i_1}'-\cdots-\nu_{i_{m}}'+1$
both of the terms
$h_r(\nu';i_1,\dots,i_{m})$
and
$h_{r-1}(\nu';i_1,\dots,i_{m})$
belong to $I_{\nu'}^\mu$.
\end{proof}

Recall $C^\mu$ denotes the algebra $C^\mu_\nu$ when $\nu$ is regular.
Let $\widetilde{C}^\mu$ denote the $\C S_n$-module obtained from
$C^\mu$ by twisting the natural action by sign.
Let $\widetilde{C}^\mu_\nu$ denote the subspace of all $S_\nu$-invariants
in $\widetilde{C}^\mu$, i.e. the space of all $S_\nu$-anti-invariants
in $C^\mu$.
The restriction of the canonical 
quotient map
$C \twoheadrightarrow C^\mu$
defines a surjective linear
map
$\widetilde{C}_\nu \twoheadrightarrow \widetilde{C}^\mu_\nu$.

\begin{Lemma}\label{grow}
The isomorphism
$\pi_*: \widetilde{C}_\nu \stackrel{\sim}{\rightarrow} C_\nu$ from
Lemma~\ref{b3} factors through the quotients to induce an isomorphism
$\bar\pi_*$
making the following diagram commute:
$$
\begin{CD}
\widetilde C_\nu &@>\pi_*>>& C_\nu\phantom{.} \\
@VVV&&@VVV\\
\widetilde C^\mu_\nu &@>>\bar\pi_*>& C^\mu_\nu
\end{CD}
$$
(Here, the vertical maps are the canonical quotient maps.)
Composing the direct sum of these isomorphisms over all $\nu$
with the isomorphism $\kappa_{\widetilde{C}^\mu}$ from Lemma~\ref{b4}, we obtain
a $\hat{\mathfrak{g}}$-module isomorphism $\bar\phi$ such that the following diagram commutes:
$$
\begin{CD}
\widehat{V}^{\otimes n} \otimes_{\C S_n} \widetilde{C} &@>\phi>>& {\textstyle\bigoplus_\nu} C_\nu\phantom{.}\\
@VVV&&@VVV\\
\widehat{V}^{\otimes n} \otimes_{\C S_n} \widetilde{C}^\mu &@>>\bar\phi>& {\textstyle\bigoplus_\nu} C^\mu_\nu.
\end{CD}
$$
(Here, $\phi$ is as in (\ref{phi}) and the vertical maps are the canonical quotients.)
\end{Lemma}

\begin{proof}
Writing $\alpha,\beta$ and $\gamma$ for the maps
induced by the canonical quotient homomorphisms,
the left hand square of the following diagram commutes:
$$
\begin{CD}
\widehat{V}^{\otimes n} \otimes_{\C S_n} \widetilde{C} &@>\kappa_{\widetilde{C}}>>& {\textstyle\bigoplus_\nu} \widetilde{C}_\nu
&@>\pi_*>>&{\textstyle\bigoplus_\nu C_\nu}\phantom{.}\\
@V\alpha VV&&@V\beta VV&&@VV\gamma V\\
\widehat{V}^{\otimes n} \otimes_{\C S_n} \widetilde{C}^\mu &@>>\kappa_{\widetilde{C}^\mu}>& {\textstyle\bigoplus_\nu} \widetilde C^\mu_\nu&@>>\bar\pi_*>&{\textstyle \bigoplus_\nu} C^\mu_\nu.\\
\end{CD}
$$
We have not yet defined the map $\bar\pi_*$ appearing in the right hand square.

Let us first 
show that the maps $\alpha$ and $\gamma\circ\pi_* \circ \kappa_{\widetilde{C}}$ have the same kernels.
Both of these maps are $\hat{\mathfrak{g}}$-module homomorphisms.
Moreover it is quite obvious that the restrictions of these maps to a regular weight space
have the same kernels.
Since any polynomial $\hat{\mathfrak{g}}$-module of degree $n$ is generated by 
any one of its regular weight spaces, this implies that the kernels are equal everywhere.
Using the commutativity of the left hand square, we deduce that the maps
$\gamma \circ \pi_*$ and $\beta$ also have the same kernels. Since they are both surjective,
this means that $\pi_*$ factors through the quotients to induce an isomorphism
$\bar\pi_*:\bigoplus_\nu \widetilde{C}^\mu_\nu \rightarrow \bigoplus_\nu C^\mu_\nu$
making the right hand square commute.
The rest of the lemma follows 
since $\phi = \pi_* \circ \kappa_{\widetilde{C}}$ and
$\bar\phi = \bar\pi_* \circ \kappa_{\widetilde{C}^\mu}$.
\end{proof}

Now we can identify the $\hat{\mathfrak{g}}$-module
$\bigoplus_\nu C^\mu_\nu$ explicitly.
This result is a natural extension of Theorem~\ref{useful},
replacing $S_n$ with $\hat{\mathfrak{g}}$.

\begin{Theorem}\label{dthm}
Let $\mu$ be a composition of $n$ with transpose partition $\lambda$.
\begin{itemize}
\item[(i)] As a $\hat{\mathfrak{g}}$-module,
$\bigoplus_\nu C^\mu_\nu$ is isomorphic to 
$\bigwedge^\mu(\widehat{V})$.
\item[(ii)] Assuming $C^\mu_\nu \neq 0$, i.e. $\lambda \geq \nu^+$,
the top graded component of $C^\mu_\nu$ is in degree
\begin{equation}
d^\mu_\nu := \sum_{i \geq 1} \lambda_i(\lambda_i-1)
- \sum_{i \in \Z} \nu_i(\nu_i-1)
\end{equation}
(which is twice the dimension of the Spaltenstein variety $F^\mu_\nu$).
\item[(iii)] As a $\hat{\mathfrak{g}}$-module,
the direct sum 
$\bigoplus_\nu C^\mu_\nu(d^\mu_\nu)$ 
of the
top graded components of all $C^\mu_\nu$
is isomorphic to 
$P^\lambda(\widehat{V})$.
\item[(iv)]
Given a non-zero vector
$x \in C^\mu_\nu$,
a regular composition $\omega$ and a composition
$\gamma$ with $\gamma^+ = \lambda$,
there exist 
operators $u,v \in U(\hat{\mathfrak{g}})$ 
and $y \in C^\mu_\omega$
such that $vx \in C^\mu_\omega$, 
$y(vx) \in C^\mu_\omega(d_\omega^\mu)$
and $u(y(vx))$ is the identity element of the one-dimensional algebra
$C^\mu_\gamma$.
\end{itemize}
\end{Theorem}

\begin{proof}
Consider
the $\hat{\mathfrak{g}}$-module isomorphism
$\bar\phi
:\widehat{V}^{\otimes n} \otimes_{\C S_n} \widetilde{C}^\mu
\rightarrow\bigoplus_\nu C^\mu_\nu$ 
from Lemma~\ref{grow}.
Combined with Theorem~\ref{useful}(i) and (\ref{hand1}), it implies that 
$\bigoplus_\nu C^\mu_\nu$ is isomorphic to $\bigwedge^\mu(\widehat{V})$
as a $\hat{\mathfrak{g}}$-module, giving (i).
Declaring that 
$\widehat{V}^{\otimes n}$ is concentrated in degree $0$,
the natural grading 
on $\widetilde{C}^\mu$
extends to a grading on
$\widehat{V}^{\otimes n} \otimes_{\C S_n} \widetilde C^\mu$,
and the action of $\hat{\mathfrak{g}}$ preserves this grading.
By Theorem~\ref{useful}(ii)
the top graded component of $\widetilde C^\mu$ is in degree
$\sum_{i \geq 1} \lambda_i (\lambda_i-1)$,
and 
by Theorem~\ref{useful}(iii) it is isomorphic to $S^\lambda$ as a $\C S_n$-module.
Hence the top graded component of
$\widehat{V}^{\otimes n} \otimes_{\C S_n} \widetilde C^\mu$ is also 
in degree $\sum_{i \geq 1} \lambda_i (\lambda_i-1)$
and by (\ref{hand2}) it is isomorphic to $P^\lambda(\widehat{V})$ as a
$\hat{\mathfrak{g}}$-module.
If $C^\mu_\nu \neq 0$, which means $\lambda \geq \nu^+$ by Lemma~\ref{van},
the $\nu$-weight space of $P^\lambda(\widehat{V})$ is non-zero, so this is
also the degree of the top graded component of the $\nu$-weight space of 
$\widehat{V}^{\otimes n} \otimes_{\C S_n} \widetilde C^\mu$.
Since the restriction of the isomorphism $\bar\phi$ to $\nu$-weight
spaces is a homogeneous map of degree $-\sum_{i \in \Z} \nu_i(\nu_i-1)$, this
proves (ii) and (iii).

Finally, we prove (iv).
The $\omega$-weight space of every irreducible polynomial representation of degree $n$
is non-zero. So any vector in any polynomial representation can be mapped to a non-zero vector
of weight $\omega$ by applying some element of $U(\hat{\mathfrak{g}})$.
So given a non-zero vector $x \in C^\mu_\nu$, there exists $v \in U(\hat{\mathfrak{g}})$ such that
$vx$ is a non-zero element of $C^\mu_\omega$.
Then by Theorem~\ref{useful}(iv) there is an element $y \in C^\mu_\omega$
so that $y(vx)$ is a non-zero element of the top graded component of $C^\mu_\omega$.
Finally by the irreducibility of 
$\bigoplus_\nu C^\mu_\nu(d^\mu_\nu)$,
we can find $u \in U(\hat{\mathfrak{g}})$ such that $u(y(vx))$ is a non-zero
vector of $C^\mu_\gamma$.
The latter algebra is
one dimensional as the $\gamma$-weight space 
of $\bigwedge^\mu(\widehat{V})$ is one dimensional,
so we can even ensure that
$u(y(vx))=1$.
\end{proof}

\begin{Corollary}\label{dcor}
The dimension of
$C^\mu_\nu$ is 
equal to the number of column strict $\lambda$-tableaux of type $\nu$.
\end{Corollary}

\begin{proof} 
The usual monomial basis
of the $\nu$-weight space of
$\bigwedge^\mu(\widehat{V})$ is parametrized in an obvious
way 
by column strict $\lambda$-tableaux of type $\nu$.
\end{proof}

\begin{Remark}\rm
We have shown that
$\bigoplus_\nu C^\mu_\nu = \bigoplus_{r \geq 0}
\bigoplus_\nu C^\mu_\nu(d^{\mu}_\nu - 2r)$, with
each
$\bigoplus_\nu C^\mu_\nu(d^{\mu}_\nu - 2r)$ being a 
$\hat{\mathfrak{g}}$-submodule.
Letting $\kappa$ be a partition of $n$ with transpose $\tau$,
the above arguments 
combined with Remark~\ref{hr} show moreover that
\begin{equation}
\sum_{r \geq 0} 
 \left[{\textstyle \bigoplus_\nu} C^\mu_\nu(d^{\mu}_\nu - 2r) : P^\kappa(\widehat{V})\right]t^r =
K_{\tau,\mu}(t),
\end{equation}
the Kostka-Foulkes polynomial.
Hence, the Hilbert polynomial of the graded algebra $C^\mu_\nu$ is given by the formula
\begin{equation}
\sum_{r \geq 0} \dim C^\mu_\nu(r) \:t^r
=
t^{d^\mu_\nu} \sum_{(\kappa,\tau)} K_{\kappa,\nu} K_{\tau,\mu}(t^{-2}),
\end{equation}
where the sum is over all pairs $(\kappa,\tau)$ of mutually transpose 
partitions of $n$.
\end{Remark}

\begin{Remark}\label{drem}\rm
Let
$i:F^\mu \hookrightarrow F$ 
and $j:F^\mu_\nu \hookrightarrow F_\nu$
be the natural inclusions.
Recall we have identified $H^*(F, \C)$ with the coinvariant algebra $C$.
In turn, by \cite{T}, we can identify 
$H^*(F^\mu, \C)$ with $C^\mu$ so that
the pull-back homomorphism $i^*:H^*(F, \C) \rightarrow H^*(F^\mu, \C)$ 
coincides with the natural quotient map $C \twoheadrightarrow C^\mu$.
Thus, the graded vector spaces $\widetilde{C}_\nu$ and $\widetilde{C}^\mu_\nu$
are identified with the
$S_\nu$-anti-invariants in $H^*(F, \C)$
and in $H^*(F^\mu, \C)$, respectively.
In \cite[Corollary 3.4(b)]{BM},
Borho and MacPherson
proved that there exists a graded
vector space isomorphism
$\widetilde{C}^\mu_\nu \stackrel{\sim}{\rightarrow} 
H^*(F^\mu_\nu, \C)$ that is homogeneous of degree $-\sum_{i \in \Z} \nu_i(\nu_i-1)$.
In fact, there is a unique such isomorphism $\bar\pi_*$
making the following diagram commute:
\begin{equation}
\begin{CD}
\widetilde C_\nu &@>\pi_*>>& H^*(F_\nu, \C)\phantom{.} \\
@Vi^* VV&&@VVj^* V\\
\widetilde C^\mu_\nu &@>>\bar\pi_*>& H^*(F^\mu_\nu, \C).
\end{CD}
\end{equation}
This last statement,
which does not seem to follow directly from \cite{BM}, will be proved in \cite{BO}.
Comparing with the first statement
of Lemma~\ref{grow}, it implies that the algebra $C^\mu_\nu$ is canonically isomorphic to the cohomology algebra
of the Spaltenstein variety $F^\mu_\nu$.
\end{Remark}

\begin{Example}\label{oneex}\rm
Let us give one small example where everything can be worked out by hand. 
Let $\mu = \nu = (\dots,0,1,2,1,0,\dots)$.
So the nilpotent matrix $x_\mu$ is simply equal to the matrix unit
$e_{2,3}$, and
the Spaltenstein variety $F^\mu_\nu$ is the space of all
partial flags $(L \subset H)$ consisting of an $e_{2,3}$-invariant 
line $L$ and an $e_{2,3}$-invariant hyperplane $H$ in $\C^4$
such that $e_{2,3} H \subseteq L$. 
As a variety, $F^\mu_\nu$ 
is isomorphic to two copies of $\mathbb{P}^2$ glued
at a point, and
its cohomology algebra is isomorphic to the algebra
$$
\C[x,y] / \langle x^3, y^3, xy \rangle
$$
with basis $\{1, x, x^2, y, y^2\}$; both these statements make good exercises.
On the other hand, using Lemma~\ref{l2},
our algebra $C^\mu_\nu$
is the quotient of $P_\nu$ by the ideal generated by all
the elementary symmetric functions in $x_1,x_2,x_3,x_4$ of positive degree
together with the elements
$\{x_1x_4, x_1x_2x_3, x_2x_3x_4\}$.
It is straightforward now by explicitly checking relations 
to see that there is an isomorphism
$H^*(F^\mu_\nu, \C) \rightarrow C^\mu_\nu$ sending 
$x \mapsto x_1, y \mapsto x_4$.
In particular, $\{1,x_1,x_1^2,x_4,x_4^2\}$ is a basis\footnote{A similar algebraic basis for 
the algebras $C^\mu_\nu$ in general
will be constructed in \cite{BO}.}
for $C^\mu_\nu$.
\end{Example}

\section{Trace maps}

In this section, we give a quite different algebraic interpretation of the
actions of the Chevalley generators $E_i$ and $F_i$ of $\hat{\mathfrak{g}}$
on $\bigoplus_\nu C_\nu$. Fix $i \in \Z$ and 
compositions $\nu,\nu'$ as in the key situation (\ref{maps}),
and set $a := \nu_i', b := \nu_{i+1}$ and $k := \sum_{j \leq i} \nu_i$.
Recall that 
the algebra $C^{S_{\nu} \cap S_{\nu'}}$ is free both as a $C_\nu$-module
with basis $1,x_k,\dots,x_k^a$
and as a $C_{\nu'}$-module with basis $1,x_k,\dots,x_k^b$. 
From now on, we reserve the notation $C_{\nu,\nu'}$ for this algebra
when viewed as a $(C_\nu,C_{\nu'})$-bimodule,
and we introduce the new notation 
$C_{\nu',\nu}$ for the same algebra when viewed
as a $(C_{\nu'},C_\nu)$-bimodule.
Of course since all our algebras are commutative, 
this distinction is quite artificial, but it helps
to keep
track of the actions (left or right) later on.
Tensoring with these bimodules
defines exact functors
\begin{align}
C_{\nu,\nu'}\otimes_{C_{\nu'}}?:
&\:\:C_{\nu'}\Mod \rightarrow C_\nu\Mod,\\
C_{\nu',\nu}\otimes_{C_\nu}?:&\:\:C_{\nu}\Mod \rightarrow C_{\nu'}\Mod,
\end{align}
where 
$A\Mod$ denotes the category of finite dimensional left $A$-modules.
The immediate goal is to prove that these functors are both left and right
adjoint to one another in canonical ways; see also \cite[Proposition 3.5]{FKS}.

\begin{Lemma}\label{riders}
There is a unique $(C_{\nu'},C_\nu)$-bimodule isomorphism
$$
\delta:
C_{\nu',\nu}
\stackrel{\sim}{\longrightarrow}
\hom_{C_\nu}(C_{\nu,\nu'}, C_\nu)
$$
such that $(\delta(x_k^{r}))(x_k^s) = (-1)^a h_{r+s-a}(\nu;i)$
for $r,s \geq 0$.
The inverse isomorphism sends
$f\mapsto (-1)^a\sum_{r=0}^{a} (-1)^{r}e_{r}(\nu';i) f(x_k^{a-r})$.
\end{Lemma}

\begin{proof}
For $0 \leq r \leq {a}$, let $\delta_r:C_{\nu,\nu'} \rightarrow C_\nu$
be the unique $C_\nu$-module homomorphism 
such that
$$
\delta_r(x_k^{{a}-s}) =\left\{
\begin{array}{ll}
0&\text{if $s < a$,}\\
(-1)^a&\text{if $s = a$.}
\end{array}\right.
$$
The maps $\delta_0,\delta_1,\dots,\delta_{a}$ form a basis for $\hom_{C_\nu}(C_{\nu,\nu'}, C_\nu)$ as a free $C_\nu$-module.
Moreover, by (\ref{mac4}), we have that
$\delta_0 (x_k^r)  = (-1)^ah_{r-a}(\nu;i)$
for any $r \geq 0$.
Hence, writing $x_k^r \delta_0$ for the map
$x \mapsto \delta_0(x x_k^r)$,
we get that
$$
x_k^{r} \delta_0 = \sum_{s=0}^{a} \delta_{s}h_{r-s}(\nu;i),
$$
since both sides map $x_k^{a-s}$ to $h_{r-s}(\nu;i)$.
Inverting using (\ref{mac1}) we get that
$$
\delta_r = \sum_{s=0}^r 
(-1)^{r-s} (x_k^s \delta_0) e_{r-s}(\nu;i)
$$
for $0 \leq r \leq{a}$.
These equations imply that 
the maps $\delta_0, x_k \delta_0, \dots, x_k^{a} \delta_0$ also form a basis
for $\hom_{C_\nu}(C_{\nu,\nu'}, C_\nu)$ as a free $C_\nu$-module.
Moreover, since 
$\sum_{s=0}^r (-1)^{r-s} x_k^se_{r-s}(\nu;i) = (-1)^r e_r(\nu';i)$
by a special case of (\ref{brad2}), we have shown that
$$
\delta_r =
(-1)^r e_r(\nu';i) \delta_0
$$
for $0 \leq r \leq{a}$.

Now we know enough to prove the lemma.
Recalling that 
$C^{S_\nu \cap S_{\nu'}} = C_{\nu,\nu'} = C_{\nu',\nu}$,
there is a well-defined $C^{S_\nu\cap S_{\nu'}}$-module homomorphism
$$
\delta:C_{\nu',\nu}
\rightarrow
\hom_{C_\nu}(C_{\nu,\nu'}, C_\nu),
\qquad 1 \mapsto \delta_0.
$$
This is automatically a 
$(C_{\nu'},C_\nu)$-bimodule homomorphism.
The elements $1, x_k,\dots,x_k^a$ form a basis for $C_{\nu',\nu}$
as a free right $C_\nu$-module, and $\delta$ maps them to to the 
functions
 $\delta_0, x_k \delta_0,\dots,x_k^a \delta_0$ which
we showed in the previous paragraph give a basis
for $\hom_{C_\nu}(C_{\nu,\nu'},C_\nu)$ as a free right $C_\nu$-module. 
Hence, $\delta$ is an isomorphism.
Moreover, 
$$
(\delta(x_k^{r}))(x_k^s) = (x_k^{r} \delta_0)(x_k^s)
= \delta_0(x_k^{r+s}) = (-1)^ah_{r+s-a}(\nu;i),
$$
so the isomorphism 
$\delta$ just constructed is precisely the map
in the statement of the lemma.
Finally, 
take any $f \in \hom_{C_\nu}(C_{\nu,\nu'},C_\nu)$.
We have that 
$$
f = (-1)^a \sum_{r=0}^{a} \delta_{r} f(x_k^{a-r})
=
(-1)^a \sum_{r=0}^{a} (-1)^{r} e_{r}(\nu';i)\delta_0f(x_k^{a-r}).
$$
Hence, $\delta^{-1}(f) = (-1)^a\sum_{r=0}^{a}
(-1)^{r} e_{r}(\nu';i)f(x_k^{a-r})$ as claimed.
\end{proof}

\begin{Corollary}\label{elf}
For any $C_\nu$-module $M$, there is a natural
$C_{\nu'}$-module isomorphism
$$
C_{\nu',\nu} \otimes_{C_\nu} M 
\stackrel{\sim}{\longrightarrow}
\hom_{C_\nu}(C_{\nu,\nu'}, M)
$$
such that $x_k^{r} \otimes m$ for any $r \geq 0$ and $m \in M$
maps to the unique $C_\nu$-module homomorphism 
sending $x_k^s$ to $(-1)^a h_{r+s-a}(\nu;i) m$ for each $s \geq 0$.
The inverse isomorphism sends $f\mapsto
(-1)^a\sum_{r=0}^{a} (-1)^{r} e_{r}(\nu';i) \otimes f(x_k^{a-r})$.
\end{Corollary}

\begin{proof}
Let $\delta$ be the isomorphism from Lemma~\ref{riders}.
We obtain the desired natural isomorphism from the composite
$$
\begin{CD}
C_{\nu',\nu} \otimes_{C_\nu} M 
&@>\delta \otimes \id_M>>
&\hom_{C_\nu}(C_{\nu,\nu'}, C_\nu)
\otimes_{C_\nu} M
\stackrel{\sim}{\longrightarrow}
\hom_{C_\nu}(C_{\nu,\nu'}, M),
\end{CD}
$$
where the second map is the obvious natural isomorphism.
Now compute this composite map and its inverse using
the explicit descriptions of
$\delta$ and $\delta^{-1}$ from Lemma~\ref{riders}.
\end{proof}

Corollary~\ref{elf} combined with adjointness of tensor and hom
implies that $(C_{\nu,\nu'}\otimes_{C_{\nu'}} ?, C_{\nu',\nu}\otimes_{C_{\nu}}?)$ is an adjoint pair of functors.
Let us write down the unit and counit of the canonical adjunction
explicitly; see also the discussion following (\ref{ht}) below.
The unit is the natural transformation
\begin{equation}\label{epsalpha}
\iota':\Id_{C_{\nu'}\Mod} \rightarrow C_{\nu',\nu}\otimes_{C_\nu}
C_{\nu,\nu'} \otimes_{C_{\nu'}}  ?
\end{equation}
defined by
$m \mapsto (-1)^a\sum_{r=0}^{a} (-1)^{r} e_{r}(\nu';i) \otimes x_k^{a-r} 
\otimes m$ for all $C_{\nu'}$-modules $M$ and $m \in M$. 
(The prime in the notation $\iota'$
is intended to help
remember that it defines
$C_{\nu'}$-module homomorphisms.)
The counit is the natural transformation
\begin{equation}\label{etaalpha}
\eps:
 C_{\nu,\nu'}\otimes_{C_{\nu'}}
C_{\nu',\nu} \otimes_{C_\nu}  ?
\rightarrow \Id_{C_\nu\Mod}
\end{equation}
defined 
by
$x_k^r \otimes x_k^{s} \otimes m \mapsto (-1)^ah_{r+s-a}(\nu;i) m$
for
 $r,s \geq 0$, $m \in M$
and any $C_\nu$-module $M$.

Repeating the proof of Lemma~\ref{riders} with the roles
of $C_\nu$ and $C_{\nu'}$ reversed, one shows
instead that there exists a unique $(C_\nu,C_{\nu'})$-bimodule isomorphism
\begin{equation}
\delta':
C_{\nu,\nu'}
\rightarrow \hom_{C_{\nu'}}(C_{\nu',\nu}, C_{\nu'})
\end{equation}
such that $(\delta'(x_k^{r}))(x_k^s) = h_{r+s-{b}}(\nu';i+1)$
for each $r,s \geq 0$.
The inverse map sends $f \mapsto \sum_{r=0}^{b}
(-1)^r e_{r}(\nu;i+1) f(x_k^{b-r})$.
Hence, just like in Corollary~\ref{elf},
the functors $C_{\nu,\nu'} \otimes_{C_{\nu'}} ?$
and $\hom_{C_{\nu'}}(C_{\nu',\nu}, ?)$ are isomorphic.
This means that
$(C_{\nu',\nu}\otimes_{C_{\nu}}?, C_{\nu,\nu'}\otimes_{C_{\nu}}?)$ is an adjoint pair too.
This way round, 
the unit of the canonical adjunction is the natural transformation
\begin{equation}\label{epsbeta}
\iota:\Id_{C_\nu\Mod} \rightarrow C_{\nu,\nu'}\otimes_{C_{\nu'}}
C_{\nu',\nu} \otimes_{C_\nu}  ?
\end{equation}
defined by
$m \mapsto \sum_{r=0}^{b} (-1)^{r} e_{r}(\nu;i+1)\otimes x_k^{b-r} \otimes m$ for all $C_\nu$-modules $M$ and $m \in M$. 
The counit is the natural transformation
\begin{equation}\label{etabeta}
\eps':
 C_{\nu',\nu}\otimes_{C_\nu}
C_{\nu,\nu'} \otimes_{C_{\nu'}}  ?
\rightarrow \Id_{C_{\nu'}\Mod}
\end{equation}
defined 
by
$x_k^r \otimes x_k^{s} \otimes m \mapsto h_{r+s-{b}}(\nu';i+1) m$
for
 $r,s \geq 0$, $m \in M$
and any $C_{\nu'}$-module $M$.

Using the adjoint pairs
just constructed
and the general construction of \cite[$\S$3]{Btr},
we can define some natural 
{\em trace maps}
\begin{equation}\label{t}
Z(C_\nu\Mod)
\quad\qquad
Z(C_{\nu'}\Mod),
\begin{picture}(0,20)
\put(-80,4){\makebox(0,0){$\stackrel{\stackrel{\,\scriptstyle{F_i}_{\phantom{,}}}{\displaystyle\longrightarrow}}{\stackrel{\displaystyle\longleftarrow}{\scriptstyle{E_i}}}$}}
\end{picture}
\end{equation}
\vspace{0mm}

\noindent
between the centers of the module categories.
To define these maps explicitly,
take $z\in Z(C_{\nu}\Mod)$
and $z'\in Z(C_{\nu'}\Mod)$.
Then
$E_i(z')$ 
and
$F_i(z)$ are the unique elements of 
$Z(C_{\nu}\Mod)$ and  
$Z(C_{\nu'}\Mod)$, respectively,
defined by the equations
\begin{align}\label{tbae}
E_i(z') &:= 
\eps \circ
(\bid_{C_{\nu,\nu'}\otimes_{C_{\nu'}}?}) z' (\bid_{C_{\nu',\nu}\otimes_{C_{\nu}}?})
\circ \iota,\\
F_i (z) &:= 
\eps' \circ
(\bid_{C_{\nu',\nu}\otimes_{C_{\nu}}?}) z (\bid_{C_{\nu,\nu'}\otimes_{C_{\nu'}}?})
\circ \iota'.\label{tbaf}
\end{align}
Since $C_\nu$ is a commutative algebras, 
the center 
$Z(C_\nu\Mod)$ 
is canonically isomorphic
to the algebra $C_\nu$ itself,
via the map
$C_\nu \stackrel{\sim}{\rightarrow} Z(C_\nu\Mod)$ arising from multiplication.
Similarly, $Z(C_{\nu'}\Mod)$ is canonically isomorphic to
$C_{\nu'}$.
Making these identifications, 
the maps $E_i$ and $F_i$ just defined become linear maps
between $C_\nu$ and $C_{\nu'}$.

\begin{Theorem}\label{benew}
The linear maps $E_i$ and $F_i$ just defined 
are the same as the maps arising from the
Chevalley generators $E_i$ and $F_i$
of $\hat{\mathfrak{g}}$ that were computed in Theorem~\ref{be}.
\end{Theorem}

\begin{proof}
We just check this in the case of $F_i$, the argument for
$E_i$ being similar.
We need the maps $$
\iota':C_{\nu'} \rightarrow
C_{\nu',\nu} \otimes_{C_\nu} C_{\nu,\nu'},\qquad
\eps':C_{\nu',\nu}\otimes_{C_\nu} C_{\nu,\nu'} \rightarrow C_{\nu'}
$$
arising from (\ref{epsalpha}) and (\ref{etabeta})
applied to the module $M = C_{\nu'}$.
The former maps $1 \mapsto (-1)^a\sum_{r=0}^{a} (-1)^r e_{r}(\nu';i) \otimes x_k^{a-r}$.
The latter is
the $C_{\nu'}$-module homomorphism 
mapping
$x_k^r \otimes x_k^s\mapsto h_{r+s-{b}}(\nu';i+1)$
for all $r,s\geq 0$.
Now take any $z \in C_\nu$ and write it as
$z = \sum_{r = 0}^{b} z_r x_k^r$ for (unique) 
elements $z_r \in C_{\nu'}$.
By the definition (\ref{tbaf}), 
$F_i (z)$ is the image of $1 \in C_{\nu'}$ under the composite 
first of
$\iota':C_{\nu'} \rightarrow C_{\nu',\nu}\otimes_{C_\nu} 
C_{\nu,\nu'}$,
then the endomorphism 
of 
$C_{\nu',\nu} \otimes_{C_\nu} C_{\nu,\nu'}$
defined by left multiplication by $1\otimes z$, and finally
the map $\eps':
C_{\nu',\nu}\otimes_{C_\nu} 
C_{\nu,\nu'} \rightarrow C_{\nu'}$.
As $z$ belongs to $C_\nu$, left multiplication by $1 \otimes z$ is the same thing
as left multiplication by $z \otimes 1$.
Using these facts, an elementary calculation now gives that
$$
F_i (z)
=\sum_{r = 0}^{b} z_r \cdot (-1)^a\sum_{s=0}^{a} (-1)^s e_{s}(\nu';i) 
h_{r-s+a-b}(\nu';i+1).
$$
This agrees with the second description of the action of the Chevalley generator $F_i$ on $z$ from Theorem~\ref{be}(i).
\end{proof}

\section{Category $\mathcal{O}$}

At last it is time to bring category $\mathcal O$ into the picture.
The basic notation concerning $\mathcal O$
was introduced already in the introduction.
We continue to write $\bigoplus_\nu$ for the direct sum over all compositions
$\nu$ of $n$, so 
$\bigoplus_\nu \mathcal O_\nu$ is 
the sum of all the integral blocks of $\mathcal O$. 
We also fix throughout the section a composition $\mu$ of $n$ and let
$\lambda$ be the transpose partition.
We will soon 
need the following result, which was proved in \cite[Theorem 2]{cyclo}.
For the final part of (ii),
we refer the reader to 
\cite[$\S$4]{schur} where an explicit parametrization of the
irreducible modules in parabolic category $\mathcal O$
for $\mathfrak{g} = \mathfrak{gl}_n(\C)$
is explained in terms of
column strict tableaux.

\begin{Lemma}\label{onto}
Let $\nu$ be a composition of $n$.
\begin{itemize}
\item[(i)] The natural multiplication map
$m_\nu^\mu:Z(\mathfrak{g}) \rightarrow Z(\mathcal O^\mu_\nu)$ is surjective.
\item[(ii)] 
The dimension of $Z(\mathcal O^\mu_\nu)$ is equal to the number
of isomorphism classes of irreducible modules in 
$\mathcal O^\mu_\nu$, which is the same as
the number of column strict $\lambda$-tableaux
of type $\nu$.
\end{itemize}
\end{Lemma}

Let us recall some of Soergel's results from \cite{Soergel}.
Fix a composition $\nu$ of $n$. 
Let $Q_\nu$ denote the antidominant 
projective indecomposable module in $\mathcal O_\nu$, that is,
the projective cover of the irreducible module $L(\alpha)$
where $\alpha = \sum_{i=1}^n a_i \eps_i$
is the unique weight such that $a_1 \leq \cdots \leq a_n$
and exactly $\nu_i$ of the integers $a_1,\dots,a_n$ are equal to
$i$ for each $i \in \Z$.
Let $p_\nu:Z(\mathfrak{g}) \rightarrow 
\End_{\mathcal{O}_\nu}(Q_\nu)^{\op}$
be the homomorphism induced by multiplication;
the $\op$ here indicates that we are for once viewing
$Q_\nu$ as a {\em right} module over
$\End_{\mathcal{O}_\nu}(Q_\nu)^{\op}$.
Also let
$q_\nu:Z(\mathfrak{g}) \rightarrow C_\nu$
be the homomorphism sending the generator
$z_r \in Z(\mathfrak{g})$ to $e_r(x_1+a_1,\dots,x_n+a_n) \in C_\nu$
for each $r=1,\dots,n$, as in the introduction.
With this notation,
we can now formulate Soergel's fundamental theorem
\cite[Endomorphismensatz 7]{Soergel} 
for the Lie algebra $\mathfrak{g} = \mathfrak{gl}_n(\C)$
as follows: 
{\em The maps $p_\nu$ and $q_\nu$ 
are surjective and have the same kernels}.
 Hence, there is a unique
isomorphism $c_\nu$ making the following diagram commute:
\begin{equation}\label{Stt}
\begin{CD}
&\!\!Z(\mathfrak{g})\:\: \\
\\
\End_{\mathcal{O}_\nu}(Q_\nu)^{\op}
@>\sim >c_\nu> &C_\nu.
\end{CD}
\begin{picture}(0,0)
\put(-66,-1){\makebox(0,0){$\swarrow$}}
\put(-68,7){\makebox(0,0){$_{p_\nu}$}}
\put(-66,-1){\line(1,1){14}}
\put(-16,-1){\makebox(0,0){$\searrow$}}
\put(-14.5,7){\makebox(0,0){$_{q_\nu}$}}
\put(-16,-1){\line(-1,1){14}}
\end{picture}
\end{equation}
We also need the following known lemma;
see for example
\cite[Theorem 5.2(2)]{MS} where a more general result than this 
is proved (for regular blocks).
For completeness, 
we include a self-contained proof based on Lemma~\ref{onto}.

\begin{Lemma}\label{iso}
The natural multiplication 
map $f_\nu:
Z(\mathcal O_\nu) \rightarrow \End_{\mathcal O_\nu}(Q_\nu)^{\op}$
is an isomorphism. 
\end{Lemma}

\begin{proof}
As well as being the projective cover of $L(\alpha)$,
the self-dual module $Q_\nu$ is its injective hull.
Every Verma module in $\mathcal O_\nu$ has irreducible
socle isomorphic to $L(\alpha)$.
Every projective module in $\mathcal O_\nu$ has a Verma flag,
so embeds into a direct sum of copies of $Q_\nu$.
Hence every module in $\mathcal O_\nu$ is a quotient of a submodule of a direct
sum of copies of $Q_\nu$.
This means that if $z$ belongs to $\ker f_\nu$, i.e. it acts as zero on $Q_\nu$,
it also acts as zero on every module in $\mathcal O_\nu$, so in fact $z = 0$.
Hence $f_\nu$ is injective.
It is surjective because the natural multiplication map 
$m_\nu:Z(\mathfrak{g}) \rightarrow Z(\mathcal O_\nu)$ is surjective
by Lemma~\ref{onto}(i), and the surjection $p_\nu$ from (\ref{Stt}) 
factors as $p_\nu = f_\nu \circ m_\nu$.
\end{proof}

Using the isomorphism $c_\nu$ from (\ref{Stt})
we will from now on {identify} the algebra
$\End_{\mathcal{O}_\nu}(Q_\nu)^{\op}$
with $C_\nu$, making  $Q_\nu$ into a right $C_\nu$-module.
Using the isomorphism $f_\nu$ from Lemma~\ref{iso}
 we will also identify
$Z(\mathcal O_\nu)$ with
$\End_{\mathcal O_\nu}(Q_\nu)^{\op}$.
So it makes sense to write simply
$Z(\mathcal O_\nu) = \End_{\mathcal O_\nu}(Q_\nu)^{\op} = C_\nu$,
and the maps $p_\nu$ and $c_\nu$ in (\ref{Stt}) have been identified with
a surjection
$p_\nu:Z(\mathfrak{g}) \rightarrow Z(\mathcal O_\nu)$
and an isomorphism $c_\nu:Z(\mathcal O_\nu) \rightarrow C_\nu$,
as we wrote in the introduction.

We introduce {Soergel's combinatorial functor}
\begin{equation}
\V_\nu:\mathcal O_\nu \rightarrow C_\nu\Mod,
\qquad
\V_\nu := \hom_{\mathcal O_\nu}(Q_\nu, ?).
\end{equation}
By \cite[Struktursatz 9]{Soergel}, given any projective module
$P$ in $\mathcal O_\nu$ and any other module $M$, the
functor $\V_\nu$
defines an isomorphism
\begin{equation}\label{sss}
\V_\nu:
\hom_{\mathcal O_\nu}(M,P) \stackrel{\sim}{\longrightarrow}
\hom_{C_\nu}(\V_\nu M, \V_\nu P).
\end{equation}
So if we let $P_\nu$ denote a minimal projective generator
for $\mathcal O_\nu$ and set
\begin{equation}
A_\nu := \End_{\mathcal{O}_\nu}(P_\nu)^{\op},
\end{equation}
then the functor $\V_\nu$ defines an algebra isomorphism
$A_\nu \stackrel{\sim}{\rightarrow} \End_{C_\nu}(\V_\nu P_\nu)^{\op}$.
It is often convenient to 
identify $A_\nu$ with $\End_{C_\nu}(\V_\nu P_\nu)^{\op}$ in this way.
We can also identify $Q_\nu$ with a unique indecomposable 
summand of $P_\nu$, so there exists
an idempotent $e_\nu \in A_\nu$ such that $Q_\nu = P_\nu e_\nu$.
It is then the case that $C_\nu = e_\nu A_\nu e_\nu$.

Suppose now that 
we are given compositions $\nu$, $\nu'$ of $n$ and
exact functors 
$\ee:\mathcal O_{\nu'} \rightarrow \mathcal O_{\nu}$
and $\ef:\mathcal O_\nu \rightarrow \mathcal O_{\nu'}$
that commute with direct sums.
Assume in addition that we are given
a $(C_{\nu'}, C_\nu)$-bimodule $C_{\nu',\nu}$, 
a $(C_{\nu}, C_{\nu'})$-bimodule $C_{\nu,\nu'}$, and a pair of
isomorphisms of functors
\begin{align}\label{kpa}
\tau&: 
\V_{\nu} \circ \ee \stackrel{\sim}{\longrightarrow}
C_{\nu,\nu'} \otimes_{C_{\nu'}} ? \circ \V_{\nu'},\\
\tau'&: 
\V_{\nu'} \circ \ef \stackrel{\sim}{\longrightarrow}
C_{\nu',\nu} \otimes_{C_\nu} ? \circ \V_\nu.
\label{kpaa}\end{align}
Let $\Adj(\ee,\ef)$ 
and
$\Adj
(C_{\nu,\nu'}\otimes_{C_{\nu'}} ?, C_{\nu',\nu} \otimes_{C_\nu} ?)$
denote 
the sets of all (not necessarily graded) adjunctions making
$(\ee, \ef)$ and
$(C_{\nu,\nu'}\otimes_{C_{\nu'}} ?, C_{\nu',\nu} \otimes_{C_\nu} ?)$, respectively, 
into
adjoint pairs of functors.
Given an adjunction belonging to
$\Adj
(C_{\nu,\nu'}\otimes_{C_{\nu'}} ?, C_{\nu',\nu} \otimes_{C_\nu} ?)$, 
there is an induced $(C_{\nu'},C_\nu)$-bimodule
isomorphism
\begin{multline}\label{ht}
\delta: C_{\nu',\nu}
= 
\hom_{C_{\nu'}}(C_{\nu'}, C_{\nu',\nu} \otimes_{C_\nu} C_\nu)
\stackrel{\sim}{\longrightarrow}\\
\hom_{C_{\nu}}(C_{\nu,\nu'}\otimes_{C_{\nu'}} C_{\nu'}, C_{\nu})
=
\hom_{C_{\nu}}(C_{\nu,\nu'}, C_\nu),
\end{multline}
where the middle isomorphism is defined by the adjunction.
Conversely, any such $(C_{\nu'},C_\nu)$-bimodule isomorphism $\delta$
defines an isomorphism 
between the functors $C_{\nu',\nu}\otimes_{C_\nu} ?$ and
$\hom_{C_\nu}(C_{\nu,\nu'}, ?)$ as in the proof of Corollary~\ref{elf}.
Hence by adjointness of tensor and hom, $\delta$ induces an adjunction
between $C_{\nu,\nu'}\otimes_{C_{\nu'}} ?$ 
and $C_{\nu',\nu} \otimes_{C_{\nu}}?$, i.e. an element of
$\Adj
(C_{\nu,\nu'}\otimes_{C_{\nu'}} ?, C_{\nu',\nu} \otimes_{C_\nu} ?)$.
The unit $\iota$ and counit $\eps$ of this induced adjunction
are characterized as follows:
\begin{itemize}
\item
$\iota$ is the natural transformation $\Id_{C_{\nu'}\Mod}
\rightarrow C_{\nu',\nu}\otimes_{C_\nu} C_{\nu,\nu'} \otimes_{C_{\nu'}} ?$
defining 
the map $1 \mapsto \sum_{i=1}^k a_i \otimes b_i$
on the regular module $C_{\nu'}$,
where $\sum_{i=1}^k a_i \otimes b_i$ is the unique element of
$C_{\nu',\nu}\otimes_{C_\nu} C_{\nu,\nu'}$ such that
$\sum_{i=1}^k (\delta(a_i))(x) b_i = x$
for all $x \in C_{\nu,\nu'}$;
\item
$\eps$ is the natural transformation
$C_{\nu,\nu'}\otimes_{C_{\nu'}} C_{\nu',\nu} \otimes ?\rightarrow \Id_{C_{\nu}\Mod}$
defining the map
$C_{\nu,\nu'} \otimes_{C_{\nu'}} C_{\nu',\nu}
\rightarrow C_\nu,
a \otimes b \mapsto (\delta(b))(a)$
on the regular module $C_\nu$.
\end{itemize}
The two constructions just described 
give mutually inverse bijections between
$\Adj(C_{\nu,\nu'}\otimes_{C_{\nu'}} ?, C_{\nu',\nu} \otimes_{C_\nu} ?)$
and the set of all $(C_{\nu'},C_\nu)$-bimodule isomorphisms
from $C_{\nu',\nu}$ to $\hom_{C_{\nu}}(C_{\nu,\nu'}, C_\nu)$.
In similar fashion, an element of $\Adj(\ee,\ef)$
defines an $(A_{\nu'}, A_\nu)$-bimodule isomorphism
\begin{equation}
\hat\delta:
\hom_{\mathcal O_{\nu'}}(P_{\nu'}, \ef P_\nu)
\stackrel{\sim}{\longrightarrow} \hom_{\mathcal O_\nu}(\ee P_{\nu'}, P_\nu).
\end{equation}
This gives a bijection between
$\Adj(\ee, \ef)$ and the set of all such $(A_{\nu'},A_\nu)$-bimodule
isomorphisms.
The following lemma 
explaining how adjunctions between 
$\ee$ and $\ef$ induce adjunctions between
$C_{\nu,\nu'}\otimes_{C_{\nu'}}?$ and $C_{\nu',\nu} \otimes_{C_{\nu}} ?$
is a general result about quotient functors in the sense of \cite[$\S$III.1]{Gab}.
We include a sketch of the proof since we will exploit 
the explicit description of the
induced unit and counit later on.

\begin{Lemma}\label{hairbrush}
There is a well-defined map
$$
T:\Adj(\ee, \ef) \rightarrow \Adj(C_{\nu,\nu'}\otimes_{C_{\nu'}}?,
C_{\nu',\nu}\otimes_{C_\nu}?)
$$
sending the adjunction in $\Adj(\ee,\ef)$ with 
unit $\hat\iota$ and counit $\hat\eps$
to the adjunction in
$\Adj(C_{\nu,\nu'}\otimes_{C_{\nu'}}?,C_{\nu',\nu}\otimes_{C_\nu}?)$
with unit $\iota$ and counit $\eps$
determined by the property that the following diagrams commute:
$$
\begin{CD}
\V_{\nu'} &@>\bid_{\V_{\nu'}} \hat\iota>>& \V_{\nu'} \ef\ee \\
@V\iota \bid_{\V_{\nu'}} VV&&@VV \tau' \bid_\ee V\\
C_{\nu',\nu}\otimes_{C_\nu} C_{\nu,\nu'} \otimes_{C_{\nu'}} \V_{\nu'}
&@<<(\bid_{C_{\nu',\nu}\otimes ?}) \tau <&C_{\nu',\nu}\otimes_{C_\nu} \V_\nu \ee,
\end{CD}
$$
$$
\begin{CD}
\V_\nu   &@<\bid_{\V_\nu} \hat\eps <<& \V_{\nu}\ee \ef \\
@A\eps \bid_{\V_\nu} AA&&@VV \tau \bid_\ef V\\
C_{\nu,\nu'}\otimes_{C_{\nu'}} C_{\nu',\nu} \otimes_{C_{\nu}} \V_{\nu}
&@<<(\bid_{C_{\nu,\nu'} \otimes ?}) \tau'<&C_{\nu,\nu'}\otimes_{C_{\nu'}} \V_{\nu'} \ef.
\end{CD}
$$
\end{Lemma}

\begin{proof}
We first construct the map $T$.
Take an adjunction between $\ee$ and $\ef$
with unit $\hat\iota$ and counit $\hat\eps$,
i.e. an element of $\Adj(\ee,\ef)$.
It defines an isomorphism
$$
\hom_{\mathcal O_{\nu'}}(Q_{\nu'}, \ef Q_\nu)
\stackrel{\sim}{\longrightarrow}
\hom_{\mathcal O_\nu}(\ee Q_{\nu'} , Q_\nu),
\qquad
g \mapsto \hat\eps_{Q_\nu} \circ \ee g.
$$
Composing on the right with the inverse of the isomorphism
\begin{multline*}
\hom_{\mathcal O_{\nu'}}(Q_{\nu'}, \ef Q_\nu)
\stackrel{\V_{\nu'}}{\longrightarrow}
\hom_{C_{\nu'}}(\V_{\nu'} Q_{\nu'}, \V_{\nu'}(\ef Q_\nu))
\stackrel{\tau'}{\longrightarrow}\\
\hom_{C_{\nu'}}(\V_{\nu'} Q_{\nu'}, C_{\nu',\nu}\otimes_{C_\nu} \V_\nu Q_\nu)
=
\hom_{C_{\nu'}}(C_{\nu'}, C_{\nu',\nu} \otimes_{C_\nu} C_\nu)
= C_{\nu',\nu}
\end{multline*}
mapping $g \mapsto \tau_{Q_\nu}' \circ \V_{\nu'} g$
and on the left with the isomorphism
\begin{multline*}
  \hom_{\mathcal O_{\nu}}(\ee Q_{\nu'},Q_\nu)
\stackrel{\V_\nu}{\longrightarrow}
\hom_{C_{\nu}}(\V_\nu (\ee Q_{\nu'}),
\V_\nu Q_\nu)
\stackrel{\tau}{\longrightarrow}\\
 \hom_{C_{\nu}}(C_{\nu,\nu'}\otimes_{C_{\nu'}} \V_{\nu'} Q_{\nu'},
\V_\nu Q_\nu)
=
\hom_{C_{\nu}}(C_{\nu,\nu'}, C_\nu)
\end{multline*}
mapping $h \mapsto \V_\nu h \circ \tau_{Q_{\nu'}}^{-1}$,
we get a $(C_{\nu'},C_\nu)$-bimodule isomorphism
$\delta:C_{\nu',\nu}\rightarrow\hom_{C_{\nu}}(C_{\nu,\nu'},C_\nu)$.
As explained just after (\ref{ht}), 
this defines an element of $\Adj(C_{\nu,\nu'} \otimes_{C_{\nu'}} ?,
C_{\nu',\nu}\otimes_{C_{\nu}} ?)$. 

Now we need to
verify that the two diagrams in the statement of the lemma commute;
we just sketch the argument for the second one.
We claim for $g \in \hom_{\mathcal O_{\nu'}}(Q_{\nu'},
\ef Q_\nu)$ and $\eps$ defined via the second diagram
that
$$
\V_\nu (\hat\eps_{Q_\nu} \circ \ee g)
=
\eps_{\V_\nu Q_\nu} \circ \left(\id_{C_{\nu,\nu'}} \otimes (\tau'_{Q_\nu}
\circ \V_{\nu'} g)\right) \circ \tau_{Q_{\nu'}},
$$
equality in $\hom_{C_\nu}(\V_\nu \ee Q_{\nu'}, \V_\nu Q_\nu)$.
Well, by the naturality of $\tau$, we have that
$\tau_{\ef Q_\nu} \circ \V_\nu \ee g = (\id_{C_{\nu,\nu'}} \otimes \V_{\nu'} g) 
\circ \tau_{Q_{\nu'}}$.
Hence, using the commuting diagram defining $\eps$
applied to the module $Q_\nu$, we get that
\begin{align*}
\V_\nu(\hat\eps_{Q_\nu} \circ \ee g)
&=
\V_\nu \hat\eps_{Q_\nu} \circ \V_\nu \ee g=
\eps_{\V_\nu Q_\nu} \circ \id_{C_{\nu,\nu'}} \otimes \tau_{Q_\nu}'
\circ \tau_{\ef Q_\nu} \circ \V_\nu \ee g \\
&= \eps_{\V_\nu Q_\nu} \circ \id_{C_{\nu,\nu'}}\otimes \tau_{Q_\nu}'
\circ \id_{C_{\nu,\nu'}}\otimes \V_{\nu'} g \circ \tau_{Q_{\nu'}}\\
&= \eps_{\V_\nu Q_\nu} \circ \left(\id_{C_{\nu,\nu'}} \otimes (\tau'_{Q_\nu} \circ
\V_{\nu'} g)\right) \circ \tau_{Q_{\nu'}}.
\end{align*}
This proves the claim.
Now take any $f \in \hom_{C_{\nu'}}(\V_{\nu'} Q_{\nu'}, C_{\nu',\nu} \otimes_{C_{\nu}} \V_\nu Q_\nu) = C_{\nu',\nu}$.
We can write $f = \tau'_{Q_\nu} \circ \V_{\nu'} g$ for a unique
$g \in \hom_{\mathcal O_{\nu'}}(Q_{\nu'}, \ef Q_\nu)$.
The map $\delta$ defined in the previous paragraph maps $f$
to $\V_\nu(\hat\eps_{Q_\nu} \circ \ee g) \circ \tau_{Q_{\nu'}}^{-1}$.
By the claim, this is the same as $\eps_{\V_\nu Q_\nu}
\circ \left(\id_{C_{\nu,\nu'}} \otimes (\tau'_{Q_\nu} \circ \V_{\nu'} g)\right)$.
This is the image of the adjunction with counit $\eps$
under the map $\delta$ from (\ref{ht}), which is what we were trying to check.
\end{proof}

Conversely, every adjunction
between $C_{\nu,\nu'}\otimes_{C_{\nu'}}?$  and
$C_{\nu',\nu}\otimes_{C_\nu} ?$ lifts in a canonical way 
to an adjunction between $\ee$ and $\ef$, thanks to the next lemma.

\begin{Lemma}\label{lift}
There exists a map
$$
R:\Adj(C_{\nu,\nu'}\otimes_{C_{\nu'}}?,
C_{\nu',\nu}\otimes_{C_\nu}?)
 \rightarrow \Adj(\ee, \ef)
$$
such that $T \circ R =\id$, where $T$ is the map from the preceeding lemma.
\end{Lemma}

\begin{proof}
Recall that there are idempotents $e_\nu \in A_\nu$
and $e_{\nu'} \in A_{\nu'}$ such that
$Q_\nu = P_\nu e_\nu, Q_{\nu'}=P_{\nu'}e_{\nu'}$
and $C_\nu = e_\nu A_\nu e_\nu, C_{\nu'} = e_{\nu'} A_{\nu'} e_{\nu'}$.
We have that
\begin{align*}
e_{\nu'}
\hom_{\mathcal O_{\nu'}}(P_{\nu'}, \ef P_\nu) e_\nu
&= \hom_{\mathcal O_{\nu'}}(Q_{\nu'}, \ef Q_\nu),\\
e_{\nu'}
\hom_{\mathcal O_{\nu'}}(\ee P_{\nu'}, P_\nu) e_\nu
&= \hom_{\mathcal O_{\nu'}}(\ee Q_{\nu'}, Q_\nu),\\
e_{\nu'}
\hom_{C_{\nu'}}(\V_{\nu'}P_{\nu'}, C_{\nu',\nu}\otimes_{C_{\nu}} \V_{\nu}P_\nu) e_\nu
&= \hom_{C_{\nu'}}(\V_{\nu'}Q_{\nu'}, C_{\nu',\nu}\otimes_{C_{\nu}} \V_{\nu}Q_\nu),\\
e_{\nu'}
\hom_{C_{\nu}}(C_{\nu,\nu'}\otimes_{C_{\nu'}}\V_{\nu'}P_{\nu'}, \V_{\nu}P_\nu) e_\nu
&= \hom_{C_{\nu'}}(C_{\nu,\nu'}\otimes_{C_{\nu'}} \V_{\nu'}Q_{\nu'}, \V_{\nu}Q_\nu).
\end{align*}
Consider the following diagram:
$$
\begin{picture}(344,200)
\put(231,153){\makebox(0,0){$\searrow$}}
\put(219,165){\line(1,-1){11}}
\put(60,105){\makebox(0,0){$\downarrow$}}
\put(59.93,162){\line(0,-1){53}}
\put(270,74){$T'$}
\put(265,55){\makebox(0,0){$\downarrow$}}
\put(264.93,100){\line(0,-1){41}}
\put(119,30){\makebox(0,0){$\searrow$}}
\put(91,58){\line(1,-1){28}}
\put(106,180){\makebox(0,0){$
\left\{
\begin{array}{l}
\text{$(A_{\nu'},A_\nu)$-bimodule isomorphisms}\\
\hom_{\mathcal O_{\nu'}}(P_{\nu'}, \ef P_\nu)
{\rightarrow}
\hom_{\mathcal O_{\nu}}(\ee P_{\nu'}, P_\nu)
\end{array}
\right\}$}}
\put(242,125){\makebox(0,0){$\left\{
\begin{array}{l}
\text{$(A_{\nu'},A_{\nu})$-bimodule isomorphisms}\\
\hom_{C_{\nu'}}(\V_{\nu'}P_{\nu'}, C_{\nu',\nu}\otimes_{C_{\nu}} \V_\nu P_\nu)
\\\qquad{\rightarrow} \hom_{C_{\nu}}(C_{\nu,\nu'}\otimes_{C_{\nu'}} \V_{\nu'}P_{\nu'}, \V_\nu P_\nu)
\end{array}
\right\}$}}
\put(105,80){\makebox(0,0){$\left\{
\begin{array}{l}
\text{$(C_{\nu'},C_{\nu})$-bimodule isomorphisms}\\
\hom_{\mathcal O_{\nu'}}(Q_{\nu'},  \ef Q_\nu)
{\rightarrow} \hom_{\mathcal O_{\nu}}(\ee Q_{\nu'}, Q_\nu)
\end{array}
\right\}
$}}
\put(242,25){\makebox(0,0){$\left\{
\begin{array}{l}
\text{$(C_{\nu'},C_\nu)$-bimodule isomorphisms}\\
\hom_{C_{\nu'}}(\V_{\nu'}Q_{\nu'}, C_{\nu',\nu}\otimes_{C_\nu} \V_\nu Q_\nu)
\\\qquad{\rightarrow}
\hom_{C_{\nu}}(C_{\nu,\nu'}\otimes_{C_{\nu'}} \V_{\nu'} Q_{\nu'}, \V_{\nu} Q_\nu)
\end{array}
\right\}
$}}
\end{picture}
$$
where the vertical maps arise by multiplying on the left by $e_{\nu'}$
and on the right by $e_\nu$, 
and the diagonal maps are the isomorphisms defined
by composing on the left and right by isomorphisms
arising from (\ref{sss}) and (\ref{kpa})--(\ref{kpaa})
exactly like in the proof of Lemma~\ref{hairbrush}.
Using the identifications of
the spaces of adjunctions with spaces of bimodule isomorphisms,
$T$ is by definition the composite of the two left hand maps.
The diagram commutes, so $T$ is also the composite of the two right hand maps.
So we just need to observe that the map 
$T'$ possesses a right inverse $R'$.
To see that, take a bimodule isomorphism $\delta$ in the codomain 
of $T'$, that is, a bimodule isomorphism
$$
\delta:C_{\nu',\nu}\stackrel{\sim}{\longrightarrow}
\hom_{C_\nu}(C_{\nu,\nu'}, C_\nu).
$$
It defines an adjunction between the functors
$(C_{\nu,\nu'} \otimes_{C_{\nu'}} ?,  C_{\nu',\nu} \otimes_{C_{\nu}} ?)$,
from which we get a bimodule isomorphism 
$$
\hat\delta:
\hom_{C_{\nu'}}(\V_{\nu'} P_{\nu'}, C_{\nu,\nu'}\otimes_{C_{\nu'}} \V_{\nu} P_\nu)
\stackrel{\sim}{\longrightarrow}
\hom_{C_{\nu}}(C_{\nu',\nu}\otimes_{C_{\nu}} \V_{\nu'} P_{\nu'}, \V_{\nu} P_\nu).
$$
Moreover, $T'(\hat\delta) = \delta$,
so we get the desired map $R'$ by setting $R'(\delta) := \hat \delta$.
\end{proof}

There is just one more essential ingredient: the special
translation functors $\ee_i, \ef_i: \mathcal O \rightarrow \mathcal O$
which we define here following the approach of \cite[$\S$7.4]{CR};
see also
\cite[$\S$3.1]{BFK} for a special case and \cite[$\S$4.4]{BKrep}
for a detailed discussion of the combinatorics of these functors in general.
Let $V$ and $V^*$ denote the natural $\mathfrak{g}$-module
and its dual, respectively.
For $\mathfrak{g}$-modules $M$ and $N$, 
multiplication by $\Omega
:= \sum_{i,j=1}^n e_{i,j} \otimes e_{j,i} \in \mathfrak{g}
\otimes \mathfrak{g}$ 
defines a $\mathfrak{g}$-module endomorphism
of $M \otimes N$.
For $i \in \Z$,
let $\ef_i$ be 
the functor defined on a module $M$ first by tensoring
with the natural module $V$, then taking the generalized $i$-eigenspace
of the endomorphism $\Omega$.
Similarly, let $\ee_i$ be the functor defined first by tensoring with
$V^*$, then taking the generalized $-(n+i)$-eigenspace of the endomorphism
$\Omega$.
This defines functors
\begin{equation}
\ee_i, \ef_i : \textstyle\bigoplus_{\nu} \mathcal O_\nu 
\rightarrow \bigoplus_{\nu} \mathcal O_\nu
\end{equation}
for each $i \in \Z$. It is well known that these functors are both
left and right adjoint to each other, for example
one gets adjunctions 
induced by the canonical adjunctions between the functors of tensoring
with $V$ and with $V^*$.
In particular, both of the functors are exact and commute with direct sums.
It is also known that $\ee_i$ is zero on modules belonging to
$\mathcal O_{\nu'}$ if $\nu_{i+1}' = 0$, and
$\ef_i$ is zero on modules belonging to $\mathcal O_{\nu}$
if $\nu_i=0$. Moreover, given compositions $\nu, \nu'$
related to each other as in the key situation (\ref{maps}), the functor $\ee_i$ maps modules belonging
to $\mathcal O_{\nu'}$ into $\mathcal O_{\nu}$ and 
$\ef_i$ maps modules belonging to $\mathcal O_\nu$ into $\mathcal O_{\nu'}$.
Hence, $\ee_i$ and $\ef_i$ restrict to well-defined
functors
\begin{equation}\label{upstairs}
\mathcal O_\nu
\quad\qquad
\mathcal{O}_{\nu'}.
\begin{picture}(0,20)
\put(-36,4){\makebox(0,0){$\stackrel{\stackrel{\,\scriptstyle{\ef_i}_{\phantom{,}}}{\displaystyle\longrightarrow}}{\stackrel{\displaystyle\longleftarrow}{\scriptstyle{\ee_i}}}$}}
\end{picture}
\end{equation}
\vspace{0mm}

From now on,
$C_{\nu,\nu'}$ denotes
the algebra $C^{S_\nu \cap S_{\nu'}}$ viewed as a 
$(C_\nu,C_{\nu'})$-bimodule
and
$C_{\nu',\nu}$ denotes 
$C^{S_\nu \cap S_{\nu'}}$ viewed as a 
$(C_{\nu'},C_\nu)$-bimodule, like in the previous section.
The following important lemma is proved in \cite[Proposition 3.3]{FKS}.

\begin{Lemma}\label{tt}
There are isomorphisms of functors
\begin{align*}
\tau&:\V_\nu \circ \ee_i \stackrel{\sim}{\longrightarrow} C_{\nu,\nu'} \otimes_{C_{\nu'}} ? \circ \V_{\nu'},\\
\tau'&:\V_{\nu'} \circ \ef_i \stackrel{\sim}{\longrightarrow} C_{\nu',\nu} \otimes_{C_{\nu}} ? \circ \V_{\nu}.
\end{align*}
\end{Lemma}

In the previous section, we constructed an explicit
adjunction between 
$C_{\nu,\nu'} \otimes_{C_{\nu'}} ?$ and $C_{\nu',\nu}\otimes_{C_{\nu}} ?$,
i.e. an element of
$\Adj(C_{\nu,\nu'}\otimes_{C_{\nu'}}?, C_{\nu',\nu}\otimes_{C_{\nu}} ?)$,
with unit $\iota'$  and counit $\eps$
defined by (\ref{epsalpha})--(\ref{etaalpha}).
Applying Lemmas~\ref{hairbrush}--\ref{lift} with $\ee = \ee_i$ and $\ef = \ef_i$
we lift this adjunction to
an element of $\Adj(\ee_i, \ef_i)$.
Denoting the unit and counit of this lift by 
$\hat\iota'$ and $\hat\eps$, respectively,
the appropriate
analogues of the diagrams in the statement of Lemma~\ref{hairbrush} 
commute.
Similarly, this time taking $\ee = \ef_i$ and $\ef = \ee_i$,
we lift the adjunction in
$\Adj(C_{\nu',\nu} \otimes_{C_{\nu}} ?, C_{\nu,\nu'}\otimes_{C_{\nu'}} ?)$
with unit $\iota$ and counit $\eps'$ defined by 
(\ref{epsbeta})--(\ref{etabeta}) to an adjunction in
$\Adj(\ef_i,\ee_i)$, 
whose unit $\hat\iota$ and counit $\hat\eps'$ are again defined
by the appropriate analogues of the diagrams from Lemma~\ref{hairbrush}.
We remark that the adjunctions
making $(\ee_i,\ef_i)$ and $(\ef_i, \ee_i)$ into adjoint pairs that we have just defined
are definitely {\em not} in general the same as 
the adjunctions induced by the
canonical adjunctions between tensoring with $V$ and $V^*$ mentioned earlier.

Now we can repeat the definitions (\ref{tbae})--(\ref{tbaf}) in the present
setting to get induced trace maps
\begin{equation}\label{egg2}
Z(\mathcal O_\nu)
\quad\qquad
Z(\mathcal{O}_{\nu'}).
\begin{picture}(0,20)
\put(-52.5,4){\makebox(0,0){$\stackrel{\stackrel{\,\scriptstyle{F_i}_{\phantom{,}}}{\displaystyle\longrightarrow}}{\stackrel{\displaystyle\longleftarrow}{\scriptstyle{E_i}}}$}}
\end{picture}
\end{equation}
\vspace{0mm}

\noindent
Thus, $E_i$ and $F_i$ are the maps
defined on
$z' \in Z(\mathcal O_{\nu'})$ 
and
$z \in Z(\mathcal O_{\nu})$ by
\begin{align}\label{newe}
E_i(z') &:= 
\hat\eps \circ \bid_{\ee_i} \,z'\, \bid_{\ef_i} \circ \hat\iota \in Z(\mathcal O_{\nu}),\\
F_i(z) 
&:= \hat\eps' \circ \bid_{\ef_i} \,z\, \bid_{\ee_i} \circ \hat\iota' \in Z(\mathcal O_{\nu'}),
\end{align}
respectively.
Also define $D_i: 
Z(\mathcal O_{\nu}) \rightarrow
Z(\mathcal O_{\nu})$ 
to be multiplication by the scalar $\nu_i$.
Taking the direct sum of these 
linear maps over all compositions of $n$, 
interpreting $E_i$ as the zero map on $Z(\mathcal O_{\nu'})$ if $\nu'_{i+1}=0$
and $F_i$ as the zero map on $Z(\mathcal O_\nu)$ if $\nu_i = 0$,
we obtain linear maps
\begin{equation}
\label{lotr}
D_i, E_i, F_i:
\textstyle\bigoplus_\nu Z(\mathcal O_\nu) \rightarrow \bigoplus_\nu Z(\mathcal O_\nu)
\end{equation}
for each $i \in \Z$.
The following theorem shows that these
maps define actions of the generators $D_i, E_i$ and $F_i$ of
$\hat{\mathfrak{g}}$ making
$\bigoplus_\nu Z(\mathcal O_\nu)$ into a well-defined
$\hat{\mathfrak{g}}$-module.

\begin{Theorem}\label{machine}
Under the identification of $\bigoplus_\nu Z(\mathcal O_\nu)$
with 
$\bigoplus_\nu C_\nu$, 
the endomorphisms $D_i, E_i$ and $F_i$ just defined
coincide with the maps arising from the actions of 
the generators $D_i, E_i$ and $F_i$ of $\hat{\mathfrak{g}}$ 
defined just after (\ref{phi}).
\end{Theorem}

\begin{proof}
It is obvious that the $D_i$'s are equal.
So in view of Theorem~\ref{benew} we
just need to check for fixed $i \in \Z$ and $\nu,\nu'$ 
as above
that the maps $E_i$ and $F_i$ from (\ref{t})
coincide with the maps $E_i$ and $F_i$ from (\ref{egg2}).
We explain the argument just for $E_i$, 
since the other case is similar.
As $C_\nu$ is commutative, we can identify
$C_\nu$ with $Z(C_\nu\Mod)$ as before.
Consider the following commutative diagram:
$$
\begin{CD}
Z(\mathfrak{g})
&@>q_\nu >>& Z(C_\nu\Mod)\\
@V p_\nu VV &&@ VV y_\nu V\\
Z(\mathcal O_\nu) &@>>x_\nu > &\End(\V_\nu),
\end{CD}
$$
where $x_\nu$ is the map sending
a natural transformation $z \in \End(\Id_{\mathcal O_\nu})$
to the natural transformation
$\bid_{\V_\nu} z \in \End(\V_\nu)$
and
$y_\nu$ is the map sending a natural transformation
$z \in \End(\Id_{C_\nu\Mod})$
to the natural transformation $z \bid_{\V_\nu} \in \End(\V_\nu)$.
We note that $y_\nu$ is injective.
Indeed, if $z \in \ker y_\nu$, then in particular 
$z$ acts as zero on
$\End_{C_{\nu}}(\V_\nu Q_\nu) = \End_{C_\nu}(C_\nu)$, hence $z = 0$.
Moreover, the diagram commutes.
To see this, take $z \in Z(\mathfrak{g})$
and any $M \in \mathcal O_\nu$.
We need to show $x_\nu(p_\nu(z))$ and $y_\nu(q_\nu(z))$
both define the same endomorphism of
$\V_\nu M = \hom_{\mathcal O_\nu}(Q_\nu,M)$.
Well, $x_\nu(p_\nu(z))$ defines the endomomorphism
$f \mapsto \hat f$ where $\hat f(q) = z f(q)$
and $y_\nu(q_\nu(z))$ defines the endomorphism
$f \mapsto \tilde f$ where $\tilde f(q) = f(zq)$.
Since $f$ is a $\mathfrak{g}$-module homomorphism we do indeed have that
$\hat f = \tilde f$.

The facts established in the previous paragraph imply that 
$x \in Z(\mathcal O_\nu)$ is equal to $y \in Z(C_\nu\Mod)$ 
under all our identifications
if and only if $\bid_{\V_\nu}x = y\bid_{\V_\nu}$
in $\End(\V_\nu)$.
Similarly, $x \in Z(\mathcal O_{\nu'})$ 
is equal to $y \in Z(C_{\nu'}\Mod)$ 
if and only if $\bid_{\V_{\nu'}}x = y\bid_{\V_{\nu'}}$
in $\End(\V_{\nu'})$.
So take $x \in Z(\mathcal O_{\nu'})$ and 
$y \in Z(C_{\nu'}\Mod)$
such that $\bid_{\V_{\nu'}} x = y \bid_{\V_{\nu'}}$.
To complete the proof of the theorem, we need to show that 
$$
\bid_{\V_\nu} E_i(x) 
= E_i(y)\bid_{\V_\nu}.
$$
Recalling (\ref{tbae}) and (\ref{newe}), 
this is the statement that
$$
\bid_{\V_\nu}\hat\eps  \circ \bid_{\V_\nu}\bid_{\ee_i} x \bid_{\ef_i}
\circ \bid_{\V_\nu}\hat \iota
=
\eps \bid_{\V_\nu} \circ (\bid_{C_{\nu,\nu'}\otimes_{C_{\nu'}}?}) 
y (\bid_{C_{\nu',\nu}\otimes_{C_{\nu}}?}) \bid_{\V_\nu} \circ 
\iota\bid_{\V_\nu}.
$$
Recalling Lemma~\ref{tt}, naturality implies that 
\begin{align*}
\tau \circ \bid_{\V_\nu} \bid_{\ee_i} x &= (\bid_{C_{\nu,\nu'}\otimes_{C_{\nu'}}?}) 
\bid_{\V_{\nu'}} x \circ \tau,\\
\tau' \circ y \bid_{\V_{\nu'}} \bid_{\ef_i}
&=
y (\bid_{C_{\nu',\nu}\otimes_{C_{\nu}}?}) \bid_{\V_\nu} \circ \tau'.
\end{align*}
Also by the commuting squares from Lemma~\ref{hairbrush} 
that are satisfied by the special adjunctions fixed above, 
we have that
\begin{align*}
\bid_{\V_{\nu}} \hat\eps 
&= \eps \bid_{\V_\nu} \circ (\bid_{C_{\nu,\nu'}\otimes_{C_{\nu'}}?}) \tau'
\circ \tau \bid_{\ef_i},\\
\iota \bid_{\V_\nu} &= 
(\bid_{C_{\nu,\nu'}\otimes_{C_{\nu'}}?})
\tau' \circ \tau \bid_{\ef_i} \circ \bid_{\V_{\nu}} \hat\iota.
\end{align*}
Now we calculate:
\begin{align*}
\bid_{\V_{\nu}} \hat\eps &\circ \bid_{\V_\nu} 
\bid_{\ee_i} x \bid_{\ef_i}
\circ \bid_{\V_{\nu}} \hat\iota\\
&=
\eps \bid_{\V_\nu} \circ (\bid_{C_{\nu,\nu'}\otimes_{C_{\nu'}}?}) \tau'
\circ \tau \bid_{\ef_i} \circ \bid_{\V_\nu} \bid_{\ee_i} x \bid_{\ef_i}
\circ \bid_{\V_\nu} \hat\iota\\
&=
\eps \bid_{\V_\nu} \circ (\bid_{C_{\nu,\nu'}\otimes_{C_{\nu'}}?}) \tau'
\circ (\bid_{C_{\nu,\nu'}\otimes_{C_{\nu'}}?})
\bid_{\V_{\nu'}} x \bid_{\ef_i}\circ \tau \bid_{\ef_i}
\circ \bid_{\V_\nu} \hat\iota\\
&=
\eps \bid_{\V_\nu} \circ (\bid_{C_{\nu,\nu'}\otimes_{C_{\nu'}}?}) \tau'
\circ (\bid_{C_{\nu,\nu'}\otimes_{C_{\nu'}}?}) 
y \bid_{\V_{\nu'}} \bid_{\ef_i}\circ \tau \bid_{\ef_i}
\circ \bid_{\V_\nu} \hat\iota\\
&=
\eps \bid_{\V_\nu} \circ (\bid_{C_{\nu,\nu'}\otimes_{C_{\nu'}}?}) 
y (\bid_{C_{\nu',\nu}\otimes_{C_{\nu}}?}) \bid_{\V_\nu} \circ 
(\bid_{C_{\nu,\nu'}\otimes_{C_{\nu'}}?}) \tau'
\circ \tau \bid_{\ef_i}
\circ \bid_{\V_\nu} \hat\iota\\
&=
\eps \bid_{\V_\nu} \circ (\bid_{C_{\nu,\nu'}\otimes_{C_{\nu'}}?}) 
y (\bid_{C_{\nu',\nu}\otimes_{C_{\nu}}?}) \bid_{\V_\nu} \circ 
\iota \bid_{\V_\nu}.
\end{align*}
This is what we wanted.
\end{proof}

Now we have all the necessary machinery set up,
we can quite quickly prove the Main
Theorem.
Let $\mathcal O^\mu_\nu$ be the integral block of 
parabolic category $\mathcal O$
parametrized by the fixed compositions $\mu$ and $\nu$,
as in the introduction\footnote{Although not needed here, we remark that Soergel's combinatorial
functor $\V$ has been considered recently in the parabolic setting in \cite[$\S$10]{Stp} and (from a quite different point of view) in \cite[$\S$5]{schur}.
}.
Because $\mathcal O_\nu^\mu$ 
is a full subcategory of $\mathcal O_\nu$, 
restriction defines an algebra homomorphism
\begin{equation}
r^\mu_\nu: Z(\mathcal O_\nu) \rightarrow Z(\mathcal O^\mu_\nu).
\end{equation} 
Since the surjection $m^\mu_\nu$ from Lemma~\ref{onto}(i)
factors as $m^\mu_\nu = r^\mu_\nu \circ p_\nu$, this map
$r^\mu_\nu$ is surjective, i.e.
$Z(\mathcal O^\mu_\nu)$ is a quotient of $Z(\mathcal O_\nu)$.
The functors $\ee_i$ and $\ef_i$ from (\ref{upstairs}) restrict to well-defined
functors
\begin{equation}\label{downstairs}
\mathcal O^\mu_\nu
\quad\qquad
\mathcal{O}^\mu_{\nu'}.
\begin{picture}(0,20)
\put(-36,4){\makebox(0,0){$\stackrel{\stackrel{\,\scriptstyle{\ef_i}_{\phantom{,}}}{\displaystyle\longrightarrow}}{\stackrel{\displaystyle\longleftarrow}{\scriptstyle{\ee_i}}}$}}
\end{picture}
\end{equation}
\vspace{0mm}

\noindent
Working with the same adjunctions as before (viewed now
as adjunctions between the restricted functors),
we define endomorphisms
\begin{equation}\label{hobbit}
D_i, E_i, F_i:
\textstyle\bigoplus_\nu Z(\mathcal O_\nu^\mu) \rightarrow
\bigoplus_\nu Z(\mathcal O_\nu^\mu)
\end{equation}
for each $i \in \Z$
in exactly the same way as (\ref{lotr}).
It is then immediate that the map 
$\oplus_\nu r^\mu_\nu:\textstyle\bigoplus_\nu Z(\mathcal O_\nu)
\twoheadrightarrow 
\bigoplus_\nu Z(\mathcal O_\nu^\mu)$
intertwines the endomorphisms
from (\ref{lotr}) and (\ref{hobbit}).
So the latter maps define actions of the generators
$D_i, E_i$ and $F_i$
of $\hat{\mathfrak{g}}$ making
$\bigoplus_\nu Z(\mathcal O^\mu_\nu)$ into a $\hat{\mathfrak{g}}$-module.
Let $s^\mu_\nu:C_\nu \twoheadrightarrow C^\mu_\nu$
denote the canonical quotient map for each $\nu$
and recall 
that $\bigoplus_\nu C^\mu_\nu$ is a $\hat{\mathfrak{g}}$-module
described by Theorem~\ref{dthm}.

\begin{Theorem}
For each composition $\nu$ of $n$, there exists a unique algebra isomorphism
$c^\mu_\nu$
making the following diagram commute:
\begin{equation}
\begin{CD}
&&
Z(\mathcal O_\nu)
= C_\nu\:\:\:
\\
\\
Z(\mathcal O_\nu^\mu)
&@>>c^\mu_\nu>&
C_\nu^\mu.
\end{CD}
\begin{picture}(0,40)
\put(-83,2){\makebox(0,0){$\swarrow$}}
\put(-85,10){\makebox(0,0){$_{r_\nu^\mu}$}}
\put(-83,2){\line(1,1){14}}
\put(-16,2){\makebox(0,0){$\searrow$}}
\put(-15.5,10){\makebox(0,0){$_{s_\nu^\mu}$}}
\put(-16,2){\line(-1,1){14}}
\end{picture}
\end{equation}
Moreover, the map
$\oplus_\nu c^\mu_\nu:
\bigoplus_\nu Z(\mathcal O^\mu_\nu) \stackrel{\sim}{\rightarrow}
\bigoplus_\nu C^\mu_\nu$ is a $\hat{\mathfrak{g}}$-module isomorphism.
\end{Theorem}

\begin{proof}
Both of the maps 
$\oplus_\nu r^\mu_\nu$
and $\oplus_\nu s^\mu_\nu$
are surjective $\hat{\mathfrak{g}}$-module homomorphisms.
Moreover, we know that $\dim Z(\mathcal O^\mu_\nu) = \dim C^\mu_\nu$
by Corollary~\ref{dcor} and Lemma~\ref{onto}(ii).
So it suffices to check that
$\ker r^\mu_\nu \subseteq \ker s^\mu_\nu$ for each $\nu$.

We first treat the special case that
$\nu^+ = \lambda$, when there is just one column strict $\lambda$-tableau
of type $\nu$.
Let $a_1 \leq \cdots \leq a_n$ be the 
integers such that $\nu_i$ of them are equal to
$i$ for each $i \in \Z$. 
Since there is just one isomorphism class of simple modules
in the highest weight category $\mathcal O^\mu_\nu$, it is a 
semisimple category.
So $z_r \in Z(\mathfrak{g})$ acts by the scalar
$e_r(a_1,\dots,a_n)$ on every module in $\mathcal O^\mu_\nu$.
It follows easily that
$\ker r^\mu_\nu$ is 
generated by the 
elements
$p_\nu(z_r) - e_r(a_1,\dots,a_n)$ 
for all $r \geq 1$.
We therefore need to show
that $
q_\nu(z_r) - e_r(a_1,\dots,a_n)
=
e_r(x_1+a_1,\dots,x_n+a_n)-e_r(a_1,\dots,a_n)$
belongs to $\ker s^\mu_\nu$ for each $r \geq 1$.
This is clear since $C^\mu_\nu$ is a one dimensional graded algebra
and each of these elements involves only strictly positive degree
terms.

Now take an arbitrary $\nu$. Also
let $\gamma$ be any composition with $\gamma^+ = \lambda$ and let
$\omega$ be any regular composition.
Let $x \in C_\nu$ be an element that is not contained in the kernel
of $s_\nu^\mu$. We need to show that $r_\nu^\mu(x) \neq 0$.
By Theorem~\ref{dthm}(iv), we can find $y \in C_\omega$ 
and $u, v \in U(\hat{\mathfrak{g}})$ such that 
$vx \in C_\omega$, 
$u(y(vx)) \in C_\gamma$ and
$s_\gamma^\mu(u(y(vx))) \neq 0$.
Since $\gamma^+=\lambda$, 
the previous paragraph
implies
that $r_\gamma^\mu(u(y(vx)) \neq 0$.
But
$$
r_\gamma^\mu(u(y(vx)) =
u r_{\omega}^\mu(y(vx)) = 
u (r_{\omega}^\mu(y) r_{\omega}^\mu(vx))=
u (r_\omega^\mu(y) (v r_\nu^\mu(x))),
$$
so we deduce that $r_\nu^\mu(x) \neq 0$ too.
\end{proof}

This completes the proof of the Main Theorem from the introduction.

\end{document}